\documentclass[titlepage,11pt]{article}
\usepackage{latexsym,tikz}
\usepackage{amsfonts}
\usepackage{amsmath}
\usepackage{textcomp}
\usetikzlibrary{graphs}
\usepackage{float}
\newcommand\blackslug{\hbox{\hskip 1pt \vrule width 4pt height 8pt depth 1.5pt
        \hskip 1pt}}
\newcommand\bbox{\hfill \quad \blackslug \bigbreak}
\newcommand{\vare}{\varepsilon}
\def\d{\hbox{-}}
\def\cc{\hbox{-}\cdots\hbox{-}}
\def\ll{,\ldots,}

\title{A survey of $\chi$-boundedness}
\author{Alex Scott\thanks{Supported by a Leverhulme Trust Research
Fellowship}\\
Mathematical Institute, University of Oxford, Oxford OX2 6GG, UK
\\
\\
Paul Seymour\thanks{Supported by ONR grant N00014-14-1-0084 and NSF
grants DMS-1265563 and DMS-1800053.}\\
Princeton University, Princeton, NJ 08544}

\newtheorem{thm}{}[section]

\newcommand{\Proof}{\noindent{\bf Proof.}\ \ }

\begin{document}
\maketitle
\begin{abstract}
If a graph has bounded clique number and sufficiently large chromatic number, 
what can we say about its induced subgraphs? Andr\'as Gy\'arf\'as made a number of challenging conjectures about this in 
the early 1980's, which have remained open until recently; but in the last  few  years there has been 
substantial  progress. This is a survey of where we are now.  
\end{abstract}

\section{Introduction}
Let $G$ be a graph. (All graphs in this paper are finite and simple.) 
We denote the chromatic number of $G$ by $\chi(G)$,
and its clique number (the cardinality of its largest clique) by $\omega(G)$.  
If $X\subseteq V(G)$, we denote the
subgraph induced on $X$ by $G[X]$, and write $\chi(X)$ for $\chi(G[X])$ when there is no danger of ambiguity.
A {\em hole} in $G$ is an
induced cycle of length at least four, and an {\em odd} hole is one with odd length; an {\em antihole} is an induced subgraph
whose complement graph is a hole of
the complement graph $\overline{G}$ of $G$.

Certainly $\chi(G)\ge \omega(G)$, and if we are told that $\chi(G) >\omega(G)$, we can deduce something about the induced
subgraphs of $G$:

\begin{thm}
If  $\chi(G) >\omega(G)$ then some induced subgraph of $G$ is an odd hole or an odd antihole. 
\end{thm}

This is the strong perfect graph theorem~\cite{SPGT}, and it settled
a long-standing open question~\cite{berge} about perfect graphs. It is in a sense the complete answer, because if we are told nothing else
about $G$ except that $\chi(G) >\omega(G)$, we cannot deduce anything more about the induced subgraphs of $G$, because $G$
might itself be the odd hole or antihole. 

But what if we are told that $\chi(G)$ is much bigger than $\omega(G)$? More precisely, fix some bound $\kappa$
and let us consider graphs $G$ with $\omega(G)\le \kappa$. If we choose $\chi(G)$ sufficiently large, can we then
deduce more about the induced subgraphs of $G$? A theorem of this type appears in~\cite{wagon}, 
but the question was first investigated systematically by 
Andr\'as Gy\'arf\'as in~\cite{gyarfastree}, and then in a beautiful and influential paper, ``Problems from the 
world surrounding perfect graphs''~\cite{gyarfas}.  Gy\'arf\'as raised a number of interesting questions.
For instance, one of his questions was: in this situation can we strengthen the 
conclusion of the strong perfect graph theorem? Perhaps such a graph must contain an odd hole?
And indeed this is true: we proved in~\cite{oddholes} that:
\begin{thm}\label{oddholes} 
For all $\kappa\ge 0$, if $G$ is a graph with $\omega(G)\le \kappa$ and 
$\chi(G)> 2^{2^{\kappa+2}}$ then $G$ has an odd hole.
\end{thm}

It is convenient to express questions of this type in the framework of {\em ideals} and {\em $\chi$-boundedness}.
An {\em ideal} (or {\em hereditary class}) is a class of graphs closed under isomorphism and under taking induced subgraphs.  
We say that a graph is {\em $H$-free} if does not contain an induced subgraph isomorphic to $H$.  Thus the class of $H$-free 
graphs is an ideal; and every ideal ${\mathcal I}$ is defined by the set of (minimal) graphs $H$ such that $\mathcal I$ is $H$-free.

An
ideal $\mathcal{I}$ is {\em $\chi$-bounded} if there is a function $f$ such that $\chi(G)\le f(\omega(G))$
for each graph $G\in \mathcal{I}$: in this case we say that $f$ is a {\em $\chi$-binding function} for $\mathcal I$.
The ideal of all graphs is not $\chi$-bounded, because as is well known (and we will see in the next section), 
there are
triangle-free graphs with arbitrarily large chromatic number; but what about subideals? Not all of them are $\chi$-bounded
(for instance, the ideal of all triangle-free graphs is not $\chi$-bounded) but some are. For instance, 
\ref{oddholes} says that the ideal of all graphs with no odd hole is $\chi$-bounded, with $\chi$-binding function 
$f(x)=2^{2^{x+2}}$.  Which other ideals are $\chi$-bounded?

There are a number of interesting conjectures and results, and to survey them we break them into three
classes: results on forests, results on holes, and other stuff.  We will cover topics roughly in that order, 
alongside related questions and topics.  We begin in the next section with some examples of graphs with 
small $\omega$ and large $\chi$; and at the end  of the paper there is a selection of open problems (in addition to those 
discussed elsewhere in the paper).

\section{Examples}\label{sec:examples}
Before we begin on any of the three main topics, let us give some useful graphs: different kinds of graphs, that are triangle-free
(that is, they have clique number two) and have arbitrarily large chromatic number. 
Most of them are explicit
constructions. (The {\em girth} of a graph is the minimum length of its cycles.)

\subsection*{Tutte's construction}

The first proof that triangle-free graphs of large chromatic number exist, is due to Tutte (writing as 
Blanche Descartes~\cite{tutte, tutte2}),
as follows.
Let $G_1$ be a 1-vertex graph. Inductively, having defined $G_k$, let $G_k$ have $n$ vertices say; now take a set $Y$ of $k(n-1)+1$
vertices, and for each $n$-subset $X$ of $Y$ take a copy of $G_k$ (disjoint from everything else), and join it to $X$ by a matching.
This makes $G_{k+1}$.
It follows inductively that $G_k$ is not $(k-1)$-colourable, and triangle-free (indeed, it has girth at least six).

\subsection*{Erd\H{o}s' random graph} 

Erd\H{o}s~\cite{erdos} proved that
for all $g,k\ge 1$, and all sufficiently large $n$, there is a graph $G$
with $n$ vertices and girth more than $G$, 
in which every stable set has fewer than $n/k$ vertices.
This can be shown as follows. Choose a function $p=p(n)$ such that $np\to \infty$ and $(np)^{g}=o(n)$. Now take 
a random graph $G$ on $2n$ vertices, in which every 
pair of vertices is joined
independently at random with probability $p$; then with probability tending to 1 as
$n\rightarrow\infty$, $G$ has no stable set of cardinality at least $n/k$, and has at most $n$ cycles 
of length 
at most $g$. Consequently, there is a set $X$ of $n$ vertices that intersects every cycle of length at 
most $g$, and by deleting $X$ we obtain a graph with the desired properties. (Erd\H{o}s actually made
a more efficient argument deleting edges instead of vertices.)
In particular, if we take $g=3$, the graph we obtain is triangle-free and has 
chromatic number more than $k$.

Later, explicit constructions for graphs with 
large girth and chromatic number were given by Lov\'asz~\cite{lovaszgirth}, Ne\v{s}et\v{r}il and R\"odl~\cite{nesetril}, and by 
Alon, Kostochka, Reiniger, West, and Zhu~\cite{AKRWZ}.

\subsection*{Mycielski's construction}

Mycielski~\cite{mycielski} gave the following construction. Let $G_2$ be the two-vertex complete graph, and 
inductively, having defined $G_{k}$, define $G_{k+1}$ as follows. Let $G_k$ have vertex set $\{v_1\ll v_n\}$.
Let $G_{k+1}$ have $2n+1$ vertices $a_1\ll a_n, b_1\ll b_n, c$, where
\begin{itemize}
\item for $1\le i<j\le n$, if $v_iv_j$ is an edge of $G_k$, then $a_ia_j, a_ib_j, b_ia_j$ are all edges
of $G_{k+1}$; and
\item for $1\le i\le n$, $b_ic$ is an edge of $G_{k+1}$.
\end{itemize}
Then $G_k$ is triangle-free and has chromatic number $k$. It is easy to see by induction that every triangle-free graph is an induced subgraph of $G_k$ for some $k$,
so this is not a good source of graphs with forbidden induced subgraphs.
  
\subsection*{Kneser graphs}

The following grew from a problem of Kneser~\cite{kneser}. Let $n,k$ be integers with $n>2k>0$, and let $K(n,k)$
be the graph with vertex set the set of all $k$-subsets of $\{1\ll n\}$, in which two such sets are adjacent
if they are disjoint. This graph has chromatic number $n - 2k + 2$, as was shown by Lov\'asz~\cite{lovasz}.
The graph $K(n,k)$ has short even cycles, but all its odd cycles have length at least $n/(n-2k)$, so if we take $n=2k(1-1/g)^{-1}$ (assuming appropriate
divisibility) we obtain a triangle-free graph in which all odd cycles have length at least $g$, and
with chromatic number at least $2k/g+2$.

\subsection*{Shift graphs}

Let $n,k$ be integers with $n>2k>2$, and let $G$ be the graph with vertex set all $k$-tuples
$(a_1\ll a_k)$ such that $1\le a_1< \cdots <a_k\le n$, in which 
$(a_1\ll a_k)$ and $(b_1\ll b_k)$ are adjacent if $a_{i+1} = b_{i}$ for $1\le i <k$, or vice versa.
This graph is triangle-free, and for fixed $k$ its chromatic number tends to infinity with $n$, 
as was shown by Erd\H{o}s and Hajnal~\cite{erdoshajnal}. Moreover, all its odd cycles have length at least
$2k+1$.

For $k=3$, this has a remarkable property. Colour each vertex $(a_1,a_2,a_3)$ by its middle element $a_2$; then this is a proper
colouring (although not optimal), and yet for every vertex, only two colours appear on its neighbours.

\subsection*{Zykov's construction}

Here is a construction due to Zykov~\cite{zykov}: let $G_1$ be a one-vertex graph, and inductively, 
having defined $G_k$, define $G_{k+1}$ as follows. Take the disjoint union of $G_1\ll G_k$
and for each $k$-tuple $(v_1\ll v_k)$
where $v_i\in V(G_i)$ for $1\le i\le k$, add a new vertex with neighbours $v_1\ll v_k$. This makes $G_{k+1}$. Then each $G_k$
is triangle-free and has chromatic number $k$.

This has a pretty variant due to Kierstead and Trotter~\cite{kiersteadtrotter}.
Take the graph $G_k$ just constructed, and orient its edges so that for 
each $i$,
the new vertices in $G_{i+1}$ are out-adjacent to their neighbours in the 
copies of $G_i$, forming a digraph $H$ say. Now for all $u,v$, if there is a 
directed path from $u$ to $v$ of odd length (more than one) in $H$, add another edge 
from $u$ to $v$. The digraph produced is still triangle-free, and it 
has no four-vertex induced directed path.

\subsection*{Ramsey graphs}

The Ramsey number $R(3,t)$ is the smallest $n$ such that every triangle-free graph with at least $n$ vertices has a $t$-vertex
stable set, and in order to bound it, constructions have been found of $n$-vertex triangle-free graphs without $t$-vertex stable sets, which therefore
have chromatic number at least $n/(t-1)$. For instance, Kim's proof~\cite{kim1} that $R(3,t)\ge O(t^2\log t)$ produced triangle-free graphs with
$O(t^2\log t)$ vertices and chromatic number $\Omega(t\log t)$. Kim's proof is non-constructive, but there are explicit 
constructions of triangle-free graphs with $\Omega(t^{3/2})$ vertices and with no stable set of 
cardinality $t$, and so with chromatic number $\Omega(t^{1/2})$, by Alon~\cite{alon}, and by Codenotti, Pudl\'ak and Resta~\cite{codenotti}.

\subsection*{The Burling graph}

Finally, here is a construction of Burling~\cite{burling}.
Let $G_1$ be the complete graph $K_2$, and let $T_1\subseteq V(G_1)$ with $|T_1|=1$. 
Inductively, suppose that we have defined $G_k$ and $T_k$, and $T_k$ is a stable subset of the vertex set of
$G_k$. We define $G_{k+1}$ and $T_{k+1}$ as follows.
Let $T_k = \{a_1\ll a_s\}$ say, and for $1\le i\le s$
let $N_i$ be the set of neighbours of $a_i$ in $G_k$.
Take a graph consisting of $s+1$ disjoint copies of $G_k\setminus T_k$, say $A_0\ll A_s$.
For $0\le i,j\le s$, let the isomorphism from $G_k\setminus T$ to $A_i$ map $N_j$ to $N_{ij}$.
Now add to this $3s^2$ new vertices, namely $x_{ij},y_{ij},z_{ij}$ for all $i,j$ with $1\le i,j\le s$.
Also add edges
so that $x_{ij},y_{ij}$ are both adjacent to every vertex in $N_{0,i}$,
and $x_{ij},z_{ij}$ are both adjacent to every vertex in $N_{ij}$, and $y_{ij}z_{ij}$ an edge, for $1\le i,j\le s$.
Let $G_{k+1}$ be the resulting graph, and let $T_{k+1}$ be the set
$$\{x_{ij},y_{ij}\::1\le i,j\le s\}.$$
It is easy to check that $G_k$ has no triangles, and for every colouring of $G_k$ with any number of colours,
some vertex in $T_k$
has neighbours of $k$ different colours, and in particular $\chi(G_k)\ge k+1$.

\section{The Gy\'arf\'as-Sumner conjecture}

An ideal is defined by the minimal graphs that it does not contain. 
If we fix a graph $H$, when is the ideal of all $H$-free graphs $\chi$-bounded?
If so, let us say
$H$ is {\em $\chi$-bounding}. For instance, the complete graph $K_2$ is obviously $\chi$-bounding, but $K_3$ is not.
Indeed, Erd\H{o}s' random graph of section~\ref{sec:examples} 
shows that no graph with a cycle is $\chi$-bounding; because
if $H$ has a cycle of length $g$, then all graphs of girth more than $g$ are $H$-free, and they are triangle-free
and can have arbitrarily large chromatic number.
Thus, only forests are $\chi$-bounding.  A famous conjecture, proposed independently
by Gy\'arf\'as~\cite{gyarfastree} and Sumner~\cite{sumner}, asserts the converse:

\begin{thm} \label{GYconj}
{\bf The Gy\'arf\'as-Sumner conjecture:} 
All forests are $\chi$-bounding.
\end{thm}

One can reduce this to trees, because a forest is $\chi$-bounding if and only if all its components are $\chi$-bounding.
To see this, let $H$ be the disjoint union of non-null forests $H_1,H_2$, which are both $\chi$-bounding; we prove
that all $H$-free graphs with clique number at most $\kappa$ have bounded chromatic number, by induction on $\kappa$.
Let $G$ be an $H$-free graph with $\omega(G)\le \kappa$; we may assume $G$ contains a subgraph $G_1$ isomorphic to $H_1$.
The set of vertices of $G$ with no neighbours in $V(G_1)$ is $H_2$-free and so has bounded chromatic number; and for each
$v\in V(G_1)$, the set of vertices adjacent to $v$ has bounded chromatic number by induction on $\kappa$. The claim follows.

Incidentally, what we just proved is a special case of something much more general. A forest that is not a tree is disconnected, and
the complete graph on $\kappa+1$ vertices is disconnected in the complement, so we could have applied the following curious
result (joint
with Maria Chudnovsky~\cite{pairs}):
\begin{thm}\label{pairs}
Let $H,J$ be graphs such that $H$ is the disjoint union of non-null graphs $H_1,H_2$, and $\overline{J}$ is the disjoint union
of non-null graphs $\overline{J_1},\overline{J_2}$. Then there exists $t$ such that for every $\{H,J\}$-free graph $G$,
its vertex set can be partitioned into $t$ subsets, such that for each of them, say $X$, $G[X]$ is either $H_1$-free, $H_2$-free,
$J_1$-free or $J_2$-free.
\end{thm}

Back to the point: for trees the conjecture is a lot tougher, and it remains open. 
The following weakening is known~\cite{scott}:
\begin{thm}\label{scott}
For every tree $T$, the ideal of all graphs containing no subdivision of $T$ as an induced 
subgraph is $\chi$-bounded.
\end{thm}

But the Gy\'arf\'as-Sumner conjecture itself 
is only known to be true for some simple types of tree. For paths there is a very pretty proof by Gy\'arf\'as~\cite{gyarfastree,gyarfas} 
that is worth giving here.\footnote{For paths, Gy\'arf\'as gives a proof for the
triangle-free case in theorem 8 of \cite{gyarfastree}; he also
says that ``J. Gerlits proved first (oral communication)'' the triangle-free case, and that Lov\'asz proved the
general case. But \cite{gyarfas} seems to be the first place a general proof is published, and this doesn't
cite anybody.}
\begin{thm}\label{paths}
Every path is $\chi$-bounding.
\end{thm}
\Proof (Sketch.) $N(v)$ denotes the neighbour set of a vertex $v$. 
We prove first the following lemma: If $G$ is a connected graph such that
\begin{itemize}
\item $\chi(N(v))\le c$ for every $v\in V(G)$; and
\item for some vertex $v$, there is no $\ell$-vertex induced path with first vertex $v$
\end{itemize}
then $\chi(G)$ is at most some function of $c,\ell$. This we prove by induction on $\ell$. Let $v$
be as in the second bullet. Then $\chi(N(v))\le c$; and for each component
$H$ of non-neighbours of $v$, some $u\in N(v)$ has a neighbour in $V(H)$, and so there is no $(\ell-1)$-vertex
induced path in $G[V(H)\cup \{v\}]$ with first vertex $u$. By the inductive hypothesis, $\chi(V(H)\cup \{u\})$ is bounded,
and so $\chi(H)$ is bounded, and therefore
$G$ has bounded chromatic number. This proves the lemma.

Then, with $\ell$ fixed, we prove by induction on $\kappa$ that for all $\kappa$, every graph $G$ with $\omega(G)\le \kappa$ and no $\ell$-vertex induced path has
bounded chromatic number. By the inductive hypothesis, $\chi(N(v))$ is bounded for each vertex $v$; 
and so from the lemma, each component of $G$, and hence $G$ itself, has bounded chromatic number.~\bbox

Until recently the complete list of cases when the Gy\'arf\'as-Sumner conjecture is known was the following:
\begin{itemize}
\item Stars (a {\em star} is a tree in which one 
vertex is adjacent to all the others); this is just an easy application of Ramsey's theorem.
\item Paths and brooms (see Gy\'arf\'as~\cite{gyarfastree,gyarfas}; a {\em broom} is obtained by identifying an 
endvertex of a path with the central vertex of a star).
\item Generalizing these, subdivisions of stars~\cite{scott} (that is, trees in which at most one vertex has degree more than two).
%(A {\em star subdivision}
%is a tree obtained from a star by subdividing edges, that is, a 
\item Trees of radius two (Gy\'arf\'as, Szemer\'edi and Tuza \cite{gst} for the triangle-free case, and 
Kierstead and Penrice~\cite{kierstead1} for the general case).
\item Trees that can be obtained from a tree of radius two by subdividing once every edge incident with the root (Kierstead and Zhu~\cite{kierstead2}).
\end{itemize}

\begin{figure}[H]
\centering

\begin{tikzpicture}[scale=.8,auto=left]
\tikzstyle{every node}=[inner sep=1.5pt, fill=black,circle,draw]
\node (z) at (0,0) {};
\node (a) at (-3,-2) {};
\node (b) at ( -1,-2) {};
\node (c) at ( 1,-2) {};
\node (d) at ( 3,-2) {};
\node (a1) at (-3.6, -3) {};
\node (a2) at (-3.0, -3) {};
\node (a3) at (-2.4, -3) {};
\node (b1) at (-1.6, -3) {};
\node (b2) at (-1.0, -3) {};
\node (b3) at (-0.4, -3) {};
\node (c1) at (0.4, -3) {};
\node (c2) at (1, -3) {};
\node (c3) at (1.6, -3) {};
\node (d1) at (2.4, -3) {};
\node (d2) at (3, -3) {};
\node (d3) at (3.6, -3) {};

\node (y) at (9,0) {};
\node (e) at (6,-2) {};
\node (f) at ( 8,-2) {};
\node (g) at ( 10,-2) {};
\node (h) at ( 12,-2) {};
\node (e1) at (5.4, -3) {};
\node (e2) at (6.0, -3) {};
\node (e3) at (6.6, -3) {};
\node (f1) at (7.4, -3) {};
\node (f2) at (8, -3) {};
\node (f3) at (8.6, -3) {};
\node (g1) at (9.4, -3) {};
\node (g2) at (10, -3) {};
\node (g3) at (10.6, -3) {};
\node (h1) at (11.4, -3) {};
\node (h2) at (12, -3) {};
\node (h3) at (12.6, -3) {};
\node (ee) at (7.5,-1) {};
\node (ff) at ( 8.5,-1) {};
\node (gg) at ( 9.5,-1) {};
\node (hh) at ( 10.5,-1) {};

\node (x) at (4.5,-4) {};
\node (i) at (1.5,-6) {};
\node (j) at ( 3.5,-6) {};
\node (k) at ( 5.5,-6) {};
\node (l) at ( 7.5,-6) {};
\node (i1) at (0.9, -7) {};
\node (i2) at (1.5, -7) {};
\node (i3) at (2.1, -7) {};
\node (j1) at (2.9, -7) {};
\node (j2) at (3.5, -7) {};
\node (j3) at (4.1, -7) {};
\node (k1) at (4.9, -7) {};
\node (k2) at (5.5, -7) {};
\node (k3) at (6.1, -7) {};
\node (l1) at (6.9, -7) {};
\node (l2) at (7.5, -7) {};
\node (l3) at (8.1, -7) {};
\node (ii) at (3,-5) {};
\node (jj) at ( 4,-5) {};

\foreach \from/\to in {z/a,z/b, z/c,z/d,a/a1,a/a2,a/a3,b/b1,b/b2,b/b3,c/c1,c/c2,c/c3,d/d1,d/d2,d/d3}
\draw [-] (\from) -- (\to);
\foreach \from/\to in {y/ee,y/ff,y/gg,y/hh,ee/e,ff/f,gg/g,hh/h,e/e1,e/e2,e/e3,f/f1,f/f2,f/f3,g/g1,g/g2,g/g3,h/h1,h/h2,h/h3}
\draw [-] (\from) -- (\to);
\foreach \from/\to in {x/ii,x/jj,x/k,x/l,ii/i,jj/j,i/i1,i/i2,i/i3,j/j1,j/j2,j/j3,k/k1,k/k2,k/k3,l/l1,l/l2,l/l3}
\draw [-] (\from) -- (\to);

\end{tikzpicture}

\caption{Kierstead-Penrice and Kierstead-Zhu trees, and our generalization} \label{fig:1}
\end{figure}
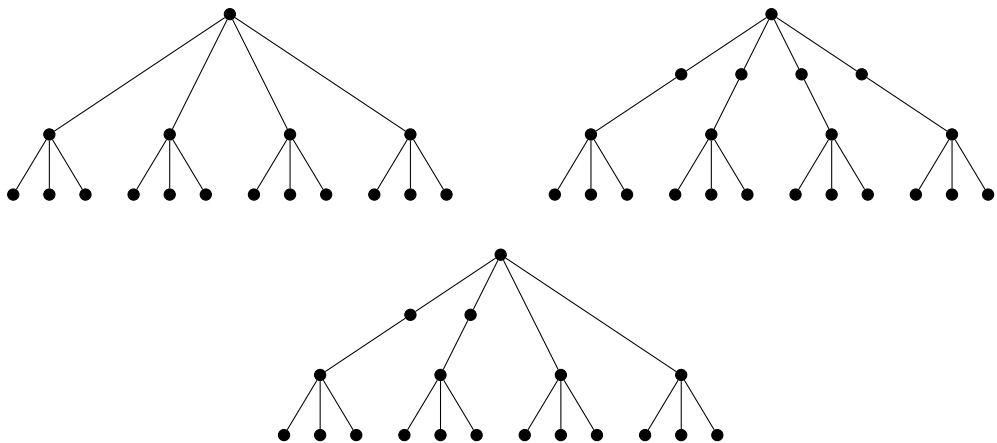

We unified the last two statements above, proving~\cite{newbrooms} that:
\begin{thm}\label{newbrooms}
Trees that can be obtained from a tree of radius two by subdividing once some edges incident with the root
are $\chi$-bounding.
\end{thm}
This uses the same proof method as~\cite{gst,kierstead1,kierstead2}, with more bells and whistles, the method of
``templates''. Here is the basic idea of that proof.
By \ref{rodltree} (see section \ref{openproblems}) and Ramsey's theorem, one can assume that
there is a $K_{n,n}$ induced subgraph, with $n$
a large constant. There might or might not be a $K_{n',n',n'}$ induced subgraph (the complete 3-partite graph), where $n'$
is a large constant (much smaller than $n$), and so on; look at the largest $k$ such that $G$ has a subgraph $G[X]$
partitioned
into $k$ stable parts, all of some specified large size and complete to one another ($k$ is bounded since the clique number
is bounded). Then we ask how the other vertices attach to $X$. Not many vertices have a large number of neighbours
in each part of $X$, from the maximality of $X$; and every other vertex $v$ with a neighbour in $X$ is a leaf of 
a large star of $G[X\cup \{v\}]$,
and these stars are useful for growing the desired tree. Delete $X$ and some appropriate
set $Y(X)$ of vertices with neighbours in $X$ (deciding which vertices to delete is the delicate part of the proof)
and do it again,
and so on; when the process stops, what remains is a set of vertices with no $K_{n,n}$ induced subgraph and with bounded 
clique number, so with bounded chromatic number by \ref{rodltree}, and we can delete it. In effect the vertex set of the graph
is partitioned into ``templates'', the various sets $X\cup Y(X)$. One can assume that each template has bounded chromatic 
number (because we can assume that for every vertex, its first neighbourhood has bounded chromatic number, by induction on
clique number; and each template $X\cup Y(X)$ is a union of boundedly many first neighbourhoods, since each $X$ has bounded
size).
So if $G$ has huge chromatic number, there must be many edges between the templates; and 
now the proof gets harder and we omit any more details.

The trees mentioned so far have the property that their vertices of degree more than two are close together, 
and indeed at this stage the 
conjecture still remained open for every tree with two vertices of degree more than two, far apart. 
But now we have several theorems proving that some such trees are $\chi$-bounding. Here are three (the first two are joint with Chudnovsky~\cite{distantstars};
the third is joint with Chudnovsky and Sophie Spirkl, and appears in Spirkl's PhD thesis~\cite{spirkl}):

\begin{thm}\label{distantstars}
The following are $\chi$-bounding:
\begin{itemize}
\item trees obtained from a star and a star subdivision by adding a path joining their centres;
\item trees obtained from a star subdivision by adding one vertex.
\item trees obtained from two disjoint paths by adding an edge between them.
\end{itemize}
\end{thm}

\begin{figure}[H]
\centering

\begin{tikzpicture}[scale=.8,auto=left]
\tikzstyle{every node}=[inner sep=1.5pt, fill=black,circle,draw]
\node (x) at (-.5,0) {};
\node (x1) at (-1.5,-2) {};
\node (x2) at ( -.5,-2) {};
\node (x3) at ( .5,-2) {};
\node (y) at ( 2.5, 0) {};
\node (y1) at ( 2, -1) {};
\node (y2) at (2.5,-1) {};
\node (y3) at (3,-1) {};

%\tikzstyle{every node}=[fill=red!20]
\foreach \from/\to in {y/y1,y/y2,y/y3}
\draw [-] (\from) -- (\to);
\foreach \from/\to in {x/x1,x/x2,x/x3,x/y}
\draw [dashed] (\from) -- (\to);

\node (z) at (5,0) {};
\node (z1) at (4,-2) {};
\node (z2) at ( 5,-2) {};
\node (z3) at ( 6,-2) {};
\node (w) at ( 8, 0) {};
\node (w1) at ( 7.5, -1) {};
\node (w2) at (9,-2) {};

%\tikzstyle{every node}=[fill=red!20]
\foreach \from/\to in {w/w1}
\draw [-] (\from) -- (\to);
\foreach \from/\to in {z/z1,z/z2,z/z3,z/w,w/w2}
\draw [dashed] (\from) -- (\to);

\node (a) at (11.5,0) {};
\node (a1) at (10,-2) {};
\node (a2) at (11,-2) {};
\node (b) at ( 12.5, 0) {};
\node (b1) at ( 13, -2) {};
\node (b2) at (14,-2) {};
\foreach \from/\to in {a/b}
\draw [-] (\from) -- (\to);
\foreach \from/\to in {a/a1,a/a2,b/b1,b/b2}
\draw [dashed] (\from) -- (\to);

\end{tikzpicture}

\caption{The trees of \ref{distantstars} (dashed lines are paths of arbitrary length)} \label{fig:3}
\end{figure}
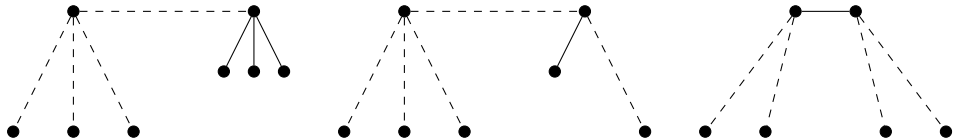

The proofs of the last three results are all by variants of the following idea. First, we work by induction on clique number,
so we may assume that the set of neighbours of each vertex has bounded chromatic number. Second, choose a vertex $v_0$, and classify
all vertices by their distance from $v_0$. (This gives us a collection of disjoint subsets $L_0,L_1,\ldots$ of $V(G)$ such that
$|L_0|=1$ and for all $i>1$, all the vertices in $L_i$ have a neighbour in $L_{i-1}$ and have no neighbours in $L_h$ for $h<i-1$; 
we call this a {\em levelling}.) In this levelling, one of the levels $L_k$ say has chromatic number at least $\chi(G)/2$, and now
we focus on it. Order the vertices in $L_{k-1}$, say $L_{k-1}=\{u_1\ll u_m\}$; and for each $i$, let $W_i$ be the set of vertices
in $L_k$ that are adjacent to $u_i$ and nonadjacent to $u_1\ll u_{i-1}$. This partitions $L_k$ into the sets $W_1,W_2\ll W_m$,
and each of the $W_i$ has bounded chromatic number (since $W_i$ has a common neighbour). But the union of all the $W_i$ has large
chromatic number, and so (for instance) there must exist some $i$ and a vertex $v\in W_i$ with many neighbours in 
$W_{i+1}\cup\cdots\cup W_m$, pairwise nonadjacent. 
Since $u_i$ is adjacent to $v$ and nonadjacent to all these neighbours, we have a little bit of a tree,
that we can combine with other parts grown elsewhere.
Incidentally, we call $(W_1\ll W_m)$ a {\em grading}; it is a surprisingly useful concept, and is used both here and in some of the results discussed in section \ref{sec:holes}.

Finally, two more results proving that some trees with far-apart vertices of degree more than two are $\chi$-bounding: 
we proved different strengthenings of the theorems of Kierstead-Penrice and of Kierstead-Zhu (though these are not
written down, so perhaps ``we believe we proved'' is more correct).
Let $T_1$ be a tree
of radius two, with root $v_1$; let $T_2$ be a star, with centre $v_2$; take the disjoint union of $T_1$ and $T_2$;
and join $v_1,v_2$ by a path, of any length. Then the tree just constructed is $\chi$-bounding. The same holds if $T_1$
is one of the Kierstead-Zhu trees. We were not able to prove that it holds when $T_1$
is one of the trees of \ref{newbrooms}.

\section{Variants of Gy\'arf\'as-Sumner: rainbow subgraphs}

Suppose $G$ is a graph with huge chromatic number and bounded clique number, and now we colour it (not necessarily optimally).
Say a subgraph is {\em rainbow} if all its vertices have different colours. Which graphs $H$ have the property that for all such $G$,
and all colourings of $G$, there is a rainbow induced copy of $H$ in $G$? Not if $H$ has a cycle (because then $H$ need not appear at
all, rainbow or not); and not if $H$ has a vertex of degree more than two, as shown by Kierstead and Trotter~\cite{kiersteadtrotter}. 
To see the latter,
let $G$ be the shift graph, described in section \ref{sec:examples}, with $k=3$. We showed earlier that $G$ can be coloured
so that the neighbours of each vertex receive only two colours, and so $G$ contains no rainbow copy of $H$ if $H$ has 
a vertex of degree more than two.

So all graphs $H$ with the desired property are paths and disjoint unions of paths; what about the converse? 
Here we have a complete solution \cite{rainbow} (see Gy\'arf\'as and S\'ark\"ozy \cite{rainbowgs} for earlier results).  As every disjoint union of paths is an induced subgraph of a larger path, it is enough to prove the following:
\begin{thm}\label{rainbowpaths}
 For every $\kappa,s\ge1$ there exists $c$ such that if $G$ is a graph with $\chi(G)\ge c$ and $\omega(G)\le \kappa$, then for
every proper colouring of $G$, there is an $s$-vertex rainbow induced path of $G$. 
\end{thm}
The proof is neat, so we sketch it. 
For a vertex $v$, let $R(v)$
be the set of vertices that are joined to $v$ by a rainbow path (or equivalently, by a rainbow induced path). 

Suppose first that for some $v$, $R(v)$ has large 
chromatic number. Enumerate the neighbours of $v$, say $u_1\ll u_m$, and for each $i$ let $U_i$ be the set of vertices 
that are reached by a rainbow induced path starting from $v$ with first edge $vu_i$. If one of these still has large chromatic number (not so large 
as before, but still large), then enumerate the neighbours of $u_i$ that are not neighbours of $v$, and repeat. 
Eventually (since we do not get a long induced rainbow
path), we find some induced path $P$ from $v$ to $u$ say, such that the set of vertices reachable by induced rainbow paths starting with $P$
has large chromatic number, but $P$ is maximal (in the sense that we cannot extend $P$ even if we somewhat reduce ``large''). 
Let $Z$ be the set of all vertices that have no neighbour in $V(P)$ except possibly $u$.
Enumerate the neighbours of $u$ in $Z$, and let $U_i$ be the set of vertices of $Z$
that are reachable by
a rainbow path that starts with $P$ and then uses the $i$th neighbour of $u$, and its further vertices all belong to $Z$ 
and are nonadjacent to $u$. For each $i$, let 
$W_i=U_i\setminus (U_1\cup\cdots\cup U_{i-1})$. This 
giving a grading $(W_1\ll W_m)$ say. Let $W=W_1\cup\cdots\cup W_m$, and direct every edge of $G[W]$ towards its end with the larger colour. 
The Gallai-Roy theorem~\cite{gallai,roy} implies
that there is a long directed path in this digraph (not induced), necessarily rainbow since it is directed;
and a theorem of Galvin, Rival and Sands~\cite{galvin} implies that
the vertex set of this directed path includes the vertex set of a large complete bipartite graph $K$, which is therefore also rainbow.
Let $w$ be the earliest vertex in $K$ (that is, belonging to $W_i$ for $i$ minimum).
It follows that $w$ has many neighbours, all of different colours, in later sets $W_j$, and not in $W_i$
(not quite; we 
have to arrange that each $W_i$ is stable, but that is easy); and then this gives a contradiction, because the rainbow path to $w$ 
could be extended to one of these neighbours, contradicting that this neighbour is not in $U_i$. 

The second case, when each set $R(v)$ has small chromatic number, is very similar
(again using Gallai-Roy and Galvin-Rival-Sands), and
we omit it.

Incidentally, there is a conjecture of Aravind (see~\cite{babu}):
\begin{thm}\label{aravind}
{\bf Conjecture:}
Let $G$ be a triangle-free graph.  Then for every colouring (not necessarily optimal) of $G$, there is a $\chi(G)$-vertex rainbow 
induced path in $G$.
%For every coloured graph $(G,\phi)$ with $\omega(G)\le 2$, there is a rainbow induced subgraph of $G$
%isomorphic to a $\chi(G)$-vertex path.
\end{thm}
This remains open, but some special cases have been proved. For instance, if we just ask for an induced path
(not necessarily rainbow), then
it holds by a theorem of Gy\'arf\'as~\cite{gyarfas}.  Or if we just ask for a rainbow path
(not necessarily induced), then it holds by the Gallai-Roy theorem~\cite{gallai, roy},
even without the bound on clique number: if we direct every edge of $G$
towards the end with higher colour, then every directed path of the digraph obtained is rainbow.
The conjecture also holds if the girth of $G$ equals its chromatic number, by a result of Babu, Basavaraju, Chandran and
Francis~\cite{babu}: in particular, if $(G,\phi)$ is a triangle-free coloured graph with $\chi(G)\ge 4$,
then some induced four-vertex path of $G$ is rainbow.

\section{Variants of Gy\'arf\'as-Sumner: orientations}

Suppose $G$ is a graph with huge chromatic number and bounded clique number, and now we direct its edges, obtaining a digraph
(which we also call $G$). Which digraphs $H$ must be present in $G$ as an induced subdigraph? Only (directed) forests, as usual;
but which ones? This question was raised by
Gy\'arf\'as~\cite{gyarfasproblem}, and he showed that the orientation
$\rightarrow \leftarrow \rightarrow$ of the four-vertex path does not have the property (to see this, take the shift graph
with $k=2$, and take its natural orientation).
Kierstead and Trotter~\cite{kiersteadtrotter} proved that the orientation $\rightarrow  \rightarrow \rightarrow $ of the same path
also does not have the property (to see this, take the variant of Zykov's construction given in section \ref{sec:examples}).

There are two other orientations of the four-vertex path, $\rightarrow \leftarrow \leftarrow$ and 
$\leftarrow \rightarrow \rightarrow$; these are equivalent up to reversing edges, so it is enough to consider the 
first. In the case of {\em acyclic} digraphs, Chv\'atal \cite{chvatal} showed that if $G$ has no induced  
$\rightarrow \leftarrow \leftarrow$ then its underlying graph is perfect; so $\chi(G)=\omega(G)$.  In the 
general case, Kierstead and R\"odl~\cite{kiersteadrodl} raised the question of whether the ideal of graphs 
that can be oriented with no induced
$\rightarrow \leftarrow \leftarrow$ is $\chi$-bounded, and Aboulker et~al.~\cite{aboulker} conjectured that it is.  We (with Chudnovsky) 
resolved the question in~\cite{orientations}:
\begin{thm}
The ideal of graphs that can be oriented with no induced
$\rightarrow \leftarrow \leftarrow$ is $\chi$-bounded.
\end{thm}
The proof is not particularly difficult, but also not very exciting. 

A much more 
challenging question is, what about stars? It is easy to see that in-directed stars and out-directed stars have the property,
but what about the star $H$ when some edges are directed in and some out? Again, this question was raised by 
Gy\'arf\'as~\cite{gyarfasproblem}.
Kierstead and R\"odl~\cite{kiersteadrodl} proved that the ideal of digraphs with no induced $H$ 
and no cyclic triangle is $\chi$-bounded; and
Aboulker et~al.~\cite{aboulker} showed that 
digraphs with no induced $H$ and no transitive triangle have bounded chromatic number 
(note that every orientation of $K_4$ has a transitive triangle, so if $G$ has no transitive triangle then $\omega(G)$ is at most 3).

We (with Chudnovsky) answered the question for stars in~\cite{orientations}:
\begin{thm}\label{distar}
If $H$ is a digraph obtained from directing the edges of a star, then the ideal of graphs that can be oriented so that no induced
subgraph is isomorphic to $H$ is $\chi$-bounded.
\end{thm}

This does not yet answer the question of which oriented trees in general have this property. Only those that do not contain
a four-vertex path directed as $\rightarrow \leftarrow \rightarrow$ or as $\rightarrow  \rightarrow \rightarrow $; but perhaps
that is the answer?  For futher discussion, see the paper of Aboulker et~al.~\cite{aboulker}.

\section{Holes}\label{sec:holes}
Along with the Gy\'arf\'as-Sumner conjecture, the best known of Gy\'arf\'as' conjectures were these three:

\begin{itemize}
\item The ideal of graphs with no odd hole is $\chi$-bounded.
\item For all $\ell\ge 0$, the ideal of graphs with no hole of length $>\ell$ is $\chi$-bounded.
\item For all  $\ell\ge 0$, the ideal of graphs with no odd hole of length $>\ell$ is $\chi$-bounded.
\end{itemize}
(The third implies the other two, of course.) These have all now been proved, in~\cite{oddholes,longholes,longoddholes}
respectively, partly in joint work with Chudnovsky and Spirkl; indeed an even stronger result is now known 
\cite{holeparity} (see \ref{holeparitya} below).
The proofs all use levellings, and focus on a level with large chromatic number.
The proof of the first was quite pretty, although too long to sketch here, and used some methods that we were
not able to use again for later papers. The second and third were both long and complicated, and both used an approach that
we used again in several other papers: to look at the maximum chromatic number of balls of bounded radius. Prove the
result assuming that every ball of some appropriate radius $r$ has bounded chromatic number; and use that as a lemma. Now
we can assume that some ball of radius $r$ has large chromatic number, as large as we like; and whenever we find an induced
 subgraph
of large chromatic number, it also must have an $r$-ball of large chromatic number. 

Let us say this more exactly. An {\em $r$-ball} in $G$ means a subset $B\subseteq V(G)$ such that every vertex in $X$
is joined by a path of $G[B]$ of length at most $r$ to some fixed vertex in $B$. If $\phi$ is a function, we say that
a graph $G$ is {\em $(r,\phi)$-controlled} 
if for every induced subgraph $H$ of $G$, $\chi(H)\le \phi(k)$ where $k$ is the maximum chromatic number of $r$-balls in $H$.
(It is helpful to assume that $\phi$ is non-decreasing, without loss of generality.)

We fix some $\kappa$ and we are trying to show that every graph $G$ with clique number at most $\kappa$ and very large chromatic number
has a hole of some desired length (let us call this a ``good'' hole).
We assume the result holds for all graphs with clique number less than $\kappa$.
Say we prove that for all graphs $G$ with clique number $\le \kappa$ and no good hole, and all $k$, if every $r$-ball has
chromatic number $\le k$, then $G$ has bounded chromatic number, say at most $\phi(k)$. That tells us that,
if $G$ is a general graph 
with clique number at most $\kappa$ and with no good hole, then $G$ is $(r,\phi)$-controlled.
This idea is very helpful and we used it many times in different situations. 

For instance, the simplest use is an idea of the first author from several years ago~\cite{scott2}, that we used to prove the theorem of~\cite{longholes}:
\begin{thm}\label{2ball}
Let $\ell\ge 0$, and let $G$ be a graph with no hole of length $>\ell$. Let every
$2$-ball in $G$ have chromatic number at most $k$. Then $G$ has chromatic number at most $4k\ell$.
\end{thm}
\Proof (Sketch.) Suppose not, choose a component with chromatic number $\chi(G)$, take a levelling of it 
with levels $L_0, L_1\ldots$, and choose $j$ such that $\chi(L_j)\ge \chi(G)/2>2k\ell$.
Let $A_1$ be a component of $G[L_j]$ with maximum chromatic number, and choose $v_0\in L_{j-1}$ with a neighbour in $A_1$. Let $A_2$
be a component with maximum chromatic number of the graph obtained from $A_1$ by deleting all neighbours of $v_0$; and 
choose a neighbour $v_1\in A_0$ of $v_0$, with a neighbour in $A_2$. Now let $A_3$ be a component with maximum chromatic
number of the graph obtained from $A_2$ by deleting all neighbours of $v_1$; and choose a neighbour $v_2\in A_1$
of $v_1$ with a neighbour in $A_3$. We continue this for $\ell$ steps. (This is Gy\'arf\'as's argument to prove \ref{paths}.)
At each stage, $\chi(A_i)$ is at least
$\chi(A_{i-1})-k$, and so, since $\chi(A_1)> 2k\ell$,
we have $\chi(A_{\ell})>k\ell$. It follows that there is some vertex in $A_{\ell}$
that has distance at least three from each of $v_0,v_1\ll v_{\ell-1}$; and this vertex has a neighbour $u\in L_{j-1}$.
Consequently, $u$ has distance at least two from each of $v_1\ll v_{\ell-1}$; and so the induced path $v_0\d v_1\cc v_{\ell-1}$
can be extended through $A_{\ell}$ to become an induced path between $v_0,u$ of length at least $\ell$. Now join $v_0,u$ by a path
via $L_0\cup\cdots\cup L_{j-2}$, and we have a hole of length $>\ell$, a contradiction. This proves \ref{2ball}.~\bbox

So far, we have wanted holes of odd length, but there are other avenues to explore. For instance, 
what about graphs with no even holes? 
Addario-Berry, Chudnovsky, Havet,  Reed and Seymour~\cite{evenholes} claimed to prove that such graphs have ``bisimplicial''
vertices, vertices whose neighbour set is the union of two cliques. Unfortunately there is a mistake in the paper and
it has been withdrawn~\cite{bisimplerratum}, but a paper~\cite{bisimplicial2} by Chudnovsky and the second author (currently being refereed)
gives a different proof. If correct, the result implies that the ideal of graphs with no even holes is $\chi$-bounded, and indeed:
\begin{thm}\label{evenholefree}
If a graph has no even hole then its chromatic number is at  most twice its clique number.
\end{thm}
The length of holes modulo 3 is also interesting. 
For instance, if $G$ is a cycle, the number of stable sets in $G$ of even cardinality              
minus the number of those of odd cardinality is $\pm 2$ if the cycle has length a multiple of three, and 
$0$ or $\pm 1$ otherwise.
Kalai and Meshulam~\cite{kalai} conjectured 
the following, which was recently proved in \cite{kmconj}:

\begin{thm}\label{stablecomplex}
If $G$ has no induced cycle of length divisible by three (and so, 
no triangles) then the number of stable sets of odd cardinality and the number of even cardinality differ by at most one.
\end{thm}

Kalai and Meshulam also conjectured that the ideal of graphs with  no induced cycle of length divisible by three 
does not contain graphs of arbitrarily large chromatic number.  This was proved in a breakthrough paper of Bonamy, 
Charbit and Thomass\'e~\cite{bonamy}:

\begin{thm}\label{bct}
Every graph with sufficiently large chromatic number contains either a triangle or a hole of
length $0$ modulo $3$.
\end{thm}
(It might be true that all these graphs are three-colourable, 
and that remains open.)  Note that \ref{bct} implies that if $\mathcal G$ is the ideal of graphs $G$ for which, 
in every induced subgraph, the number of stable sets of odd cardinality and the number of even cardinality 
differ by at most one, then the graphs in $\mathcal G$ have bounded chromatic number.
Motivated by topological considerations,
Kalai and Meshulam~\cite{kalai} made several other beautiful conjectures
connecting graph-theoretic properties with homological ones, some of which we discuss later.

More recently we proved a much stronger theorem that contains all three of the Gy\'arf\'as conjectures and \ref{bct}. We have the following~\cite{holeparity}:

\begin{thm}\label{holeparitya}
For all integers $k\ge 0$ and $\ell\ge 1$, the ideal of all graphs with no hole of length $k$ modulo $\ell$ is $\chi$-bounded.
\end{thm}

In fact, we could prove something even more general:

\begin{thm}\label{multiholeparity}
Let $n\ge0$ be an integer, and for $1 \le i \le n$ let $p_i \ge 0$ and $q_i \ge 1$ be integers. Let $\mathcal{C}$ 
be the ideal of all graphs that do not contain $n$ pairwise anticomplete holes $H_1,\dots, H_n$ where $H_i$ has 
length $p_i$ modulo $q_i$ for $1 \le i \le n$. Then $\mathcal{C}$ is $\chi$-bounded. 
\end{thm}

As an application, we used this in~\cite{holeparity}
to prove two further conjectures of Kalai and Meshulam~\cite{kalai}.  For a graph $G$, we write $I(G)$ for the 
independence complex of $G$ (that is, the collection of stable sets in $G$); the {\em Euler characteristic} of $I(G)$ 
is the number of stable sets of even cardinality minus the number of odd cardinality.  
The following, one of the Kalai-Meshulam conjectures, is proved in~\cite{holeparity}:

\begin{thm}\label{holeparityb}
For every integer $k \ge 0$ there exists $c$ such that the following holds. Let $G$ be a graph such that 
the Euler characteristic of every induced subgraph of $G$ has modulus at most $k$. 
Then $G$ has chromatic number at most $c$.
\end{thm}

The other conjecture of Kalai and Meshulam concerns the ``total Betti number'' of a graph, that is, the 
sum of the Betti numbers:

\begin{thm}\label{holeparityc}
For every integer $k \ge 0$ there exists $c$ such that the following holds. If the total Betti number of 
every induced subgraph of $G$ is at most $k$, then $\chi(G) \le c$.
\end{thm}

In both \ref{holeparityb} and \ref{holeparityc}, the key point was to look for sets of pairwise anticomplete 
holes of length divisible by three.

Despite the generality of \ref{holeparitya}, its proof is reasonable.\footnote{One of the referees said: ``Reading the proofs was even 
occasionally an enjoyable activity'' -- we are very proud of that!}
Again, we look at
the chromatic number of $r$-balls, and here $r=8$ turns out to be important. If all $8$-balls have bounded chromatic number,
we obtain a bound on $\chi(G)$ by an argument using gradings, chaining together in a cycle 
many pairs of paths whose lengths differ by one modulo $\ell$. Next we look at 7-balls; if all 7-balls have bounded chromatic
number, but some 8-ball has huge chromatic number, again we use gradings to win (with more-or-less the same argument, except 
now we can do a little better; the pairs of paths have lengths differing by one, not just by one modulo $\ell$); and so on, until
we have proved a bound on $\chi(G)$ if all 2-balls have bounded chromatic number. Now we are in much the same situation as
when we were using \ref{2ball} to prove the theorem of~\cite{longholes}, and the proof is completed the same way.  The proof of \ref{multiholeparity} uses the same ideas, with a little more work.

\ref{holeparitya} is very strong, but there are even stronger statements that might be true. We cannot hope to specify the
lengths of holes exactly, because there are graphs with large chromatic number and arbitrarily large girth, but perhaps a graph 
with huge chromatic number and bounded clique number must contain many holes with consecutive lengths. 
More exactly:
\begin{thm}\label{holeseqconj}{\bf Conjecture:}
For every integer $\ell\ge 0$, 
the ideal of graphs $G$ such that $G$ does not contain holes of $\ell$ consecutive lengths is $\chi$-bounded.
\end{thm}
This remains open,
but it is true in the triangle-free case. We proved~\cite{holeseq}:
\begin{thm}\label{holeseq}
For every integer $\ell\ge 0$, there exists $k$ such that if $G$ is triangle-free and $\chi(G)>k$, then $G$ has $\ell$ holes
of consecutive lengths.
\end{thm}
The proof of this was surprisingly difficult, and so far we see no way to extend it to graphs that have triangles.

\section{Subdivisions}\label{subdivide}

One way to formulate the result of~\cite{longholes} is: 
\begin{thm}\label{longholes}
For every cycle $C$, the ideal of all graphs containing no subdivision
of $C$ as an induced subgraph is $\chi$-bounded.
\end{thm}
We already mentioned the result from \cite{scott}, that:
\begin{thm}\label{scottagain}
For every tree $T$, the ideal of all graphs containing no 
subdivision of $T$ as an induced
subgraph is $\chi$-bounded. 
\end{thm}
What about other graphs, not cycles or trees?
Let us say $H$ is {\em weakly pervasive} if it has the property of these two theorems, that
the ideal of all graphs containing no 
subdivision of $H$ as an induced
subgraph is $\chi$-bounded. 
The first author~\cite{scott} made the conjecture that, in fact, all graphs are weakly pervasive.

Unfortunately this is false; Pawlik, Kozik, Krawczyk, Laso\'{n},
Micek, Trotter and Walczak~\cite{sevenpoles} showed that 
the Burling graph (described in 
section \ref{sec:examples}) is a counterexample. 
If $H$ is obtained from a graph $G$ by subdividing every edge exactly $\ell$ times, we say that $H$ is an {\em $\ell$-subdivision} 
of $G$; and similarly if every edge is subdivided at least $\ell$ times, we call it a {\em $(\ge \ell)$-subdivision}).
Let $H$ be a 1-subdivision of $K_5$; then the Burling graph contains no subdivision of $H$ as an induced subgraph, 
and so $H$ is not weakly pervasive, and this
disproves the conjecture.

It would be nice to characterize all the weakly pervasive graphs, but this is still open. A more tractable question
is, which graphs $H$ have the property that every subdivision of $H$ is weakly pervasive?
We call such a graph $H$ {\em pervasive}. (And extend the definition to multigraphs, for convenience.)
This we can come much closer to answering.

Note that \ref{longholes} above is equivalent to saying that the multigraph
with two vertices and two parallel edges is pervasive. Let us say a multigraph is a {\em banana tree}
if it can be obtained by adding parallel edges to a tree. 
In~\cite{bananatrees}, we proved the following, which contains both \ref{longholes}
and \ref{scottagain}:
\begin{thm}\label{bananatrees}
Every banana tree is pervasive.
\end{thm}

In the reverse direction, which graphs do we know not to be pervasive? We have seen that $K_5$ is not, because its 1-subdivision
is not contained as an induced subdivision in the Burling graph.
In fact there are smaller graphs with the same property.
For instance, let $H$ be obtained from $K_4$ by 
subdividing once every edge of a 
cycle of length four; then Chalopin, Esperet, Li and Ossona de Mendez proved in~\cite{chandeliers} that the Burling graph
contains no induced subdivision of $H$. 
Indeed, in that paper they figured out exactly which graphs $H$ had the property that, for every subdivision $H'$ of $H$,
the Burling graph contains an induced subdivision of $H'$. 
We need to describe their result.

Say a {\em chandelier} is a graph obtained from a tree by adding a new vertex (called the {\em pivot}) adjacent to all the leaves of 
the tree. A {\em tree of chandeliers} is a graph that can be obtained recursively, by starting with a 
chandelier $G_1$, 
and at each step identifying the pivot of some new chandelier with a vertex of the graph $G_i$ that we have already constructed, forming $G_{i+1}$.
A {\em forest of chandeliers} is a graph in which every component is a tree of chandeliers. It is proved in~\cite{chandeliers}
that:
\begin{thm}\label{chandeliers}
A graph $H$ has the property that, for every subdivision $H'$ of $H$, a (large enough) Burling graph contains an induced subdivision of $H'$,
if and only if $H$ is a forest of chandeliers.
\end{thm}

In view of
\ref{chandeliers}, every graph that is pervasive in the ideal of all graphs must be a forest of chandeliers.  This suggests:
\begin{thm}\label{chandconj}
{\bf Conjecture: }A graph is pervasive in the ideal of all graphs if and only if it is a forest of chandeliers.
\end{thm}
Conveniently, every subdivision of a forest of chandeliers is another forest of chandeliers, so to prove that they are
all pervasive, it is enough to prove that they are all weakly pervasive.
We have some evidence for this, in addition to the banana trees result above.

We recall that $\phi$ is a function, 
a graph $G$ is {\em $(r,\phi)$-controlled}
if for every induced subgraph $H$ of $G$, $\chi(H)\le \phi(k)$ where $k$ is the maximum chromatic number of $r$-balls in $H$.
Say an ideal $\mathcal{I}$ is {\em $r$-controlled} if there is a function $\phi$ such that 
every graph in the ideal is $(r,\phi)$-controlled.
We (with Chudnovsky) proved in~\cite{strings} that:
\begin{thm}\label{strings}
Every forest of chandeliers is pervasive in every $r$-controlled ideal, for all $r\ge 2$.
\end{thm}

In fact we proved more than this, we showed the next two results:
\begin{thm}\label{3control}
Let $m\ge 0$, let
$r\ge 2$, and let $\mathcal{I}$ be an $r$-controlled ideal. The ideal of all graphs in $\mathcal{I}$
that do not contain the $s$-subdivision of the complete bipartite graph $K_{m,m}$ as an induced subgraph, for $1\le s\le r+2$, is 
$2$-controlled. Consequently, for all $H$ and $r\ge 2$, every $r$-controlled ideal of $H$-subdivision-free graphs is $2$-controlled.
\end{thm}

\begin{thm}\label{2control} Let $\kappa, m\ge 0$, and let $H$ be a forest of chandeliers.
Let $\mathcal{I}$ be a $2$-controlled ideal of graphs all with clique number at most $\kappa$.
Then every graph in $\mathcal{I}$ with sufficiently large chromatic number
contains either $H$ or the 1-subdivision of $K_{m,m}$ as an induced subgraph.
\end{thm}

The first of these reduces our problem from $r$-controlled ideals to 2-controlled ideals (because if $G$ contains 
 an induced $s$-subdivision of the complete bipartite graph $K_{m,m}$ for $m$ large enough, then it contains an induced subdivision
of any fixed graph); and the second says that in the 2-controlled case we can do better than just getting an induced subdivision of $H$;
we can either get $H$ itself as an induced subgraph, or the 1-subdivision of $K_{m,m}$. (In fact we proved \ref{2control}
for graphs $H$ that are much more general than forests of chandeliers, but we omit the details here.)
Together, \ref{3control} and \ref{2control} imply that forests of chandeliers are weakly pervasive in any $r$-controlled ideal, and
so prove \ref{strings}. 

The proof of \ref{3control} was straightforward Ramsey theory, but \ref{2control} was intricate. It used a refinement of the 
methods used in~\cite{longholes}; let us sketch some of the ideas. We can assume there is a 2-ball with huge chromatic number,
say with centre $z_1$. The neighbours of $z_1$ have bounded chromatic number, by induction on clique number; so the set $B_1$ of 
vertices with distance exactly two from $z_1$ has huge chromatic number. In $G[B_1]$, there is therefore a 2-ball with huge (not quite so huge) chromatic number, say with centre $z_2$; and again we look at its second neighbours $B_2$ in $G[B_1]$. This contains a 2-ball
with huge chromatic number, and so on. Eventually we have a sequence of vertices $z_1, z_2\ll z_m$,
pairwise nonadjacent, and a set of neighbours $A_i$ of $z_i$ for $1\le i\le m$, pairwise disjoint, and a set $C$ with large 
chromatic number, disjoint from everything else, such that each $A_i$ covers $C$, and each $A_i$ covers $A_j\cup \{z_j\}$
for $j>i$, and $z_i$ has no neighbours in $A_j\cup \{z_j\}$ for $j>i$. Now $z_j$ has neighbours in $A_i$. Perhaps the vertices
in $A_i$ nonadjacent to $z_j$ still cover a large part of $C$, or perhaps not; and by Ramsey's theorem, we can assume the same
happens for all pairs $i<j$. In one case we can arrange to have the same structure as before except that now $z_j$ is anticomplete to 
$A_i$ for $i\ne j$; and this case is fairly straightforward, and we obtain a 1-subdivision of $K_{m,m}$. In the other case,
we can arrange that every vertex in $A_i$ is adjacent to $z_j$ for $j>i$ (at least in the triangle-free case; when there are
triangles, we need to replace the $z_i$'s by cliques and look at the common neighbours of these cliques). After a lot of intricate
analysis, we find a surprisingly rich supply of induced subgraphs, and in particular we can find $H$ itself as an induced subgraph.

A {\em string graph} is the intersection graph of a set of curves in the plane. String graphs are particularly nice in this context,
for three reasons:
\begin{itemize}
\item Pawlik, Kozik, Krawczyk, Laso\'{n}, Micek, Trotter and Walczak~\cite{sevenpoles} showed that the Burling graph is a string graph;
\item the ideal of all string graphs is 2-controlled (see \ref{stringcontrol}); and
\item no string graph contains a $(\ge 1)$-subdivision of $K_{3,3}$ as an induced subgraph.
\end{itemize}
Consequently, we have
\begin{thm}\label{stringgraphs}
For every forest of chandeliers $H$ and all $\kappa\ge 0$, every string graph $G$ with $\omega(G)\le \kappa$ and $\chi(G)$
sufficiently large contains $H$ as an induced subgraph.
\end{thm}
This is much more than saying that $H$ is pervasive in the ideal of string graphs; we get not only
an induced subdivision of $H$, but $H$ itself.

What happens for ideals that are not $r$-controlled? 
%Say an ideal is ``$r$-bounded'' if there exists $c$ such that all graphs in 
%the ideal with radius at most $r$ 
%have chromatic number at most $c$. 
%If an ideal is not $r$-controlled, then for some $c$, the ideal contains graphs with 
%arbitrarily large chromatic 
%number whose $r$-balls all have chromatic number at most $c$, which means that some $r$-bounded sub-ideal contains graphs with 
%unbounded chromatic number; so 
%an ideal is $r$-controlled if and only if every $r$-bounded subideal has bounded chromatic number.
%Hence we can assume that for any convenient, fixed value of $r$, 
%we have an ideal that is $r$-bounded and has unbounded chromatic number.
%To show that a forest of chandeliers $H$ is pervasive, we need to show there is a value of $r$ where we can handle
%$r$-bounded ideals, and that motivates the following definition.
Let us say a multigraph $H$ is {\em widespread} if for every subdivision $H'$
of $H$, there exists $r\ge 0$ such that
for every ideal of $H'$-subdivision-free graphs is $r$-controlled.
Being widespread is roughly complementary to being pervasive in $r$-controlled ideals;
a graph is both pervasive in all $r$-controlled ideals for all $r$, and widespread, if and only if it is pervasive
in the ideal of all graphs. 

\begin{figure}[H]
\centering

\begin{tikzpicture}[scale=.5,auto=left]
\tikzstyle{every node}=[inner sep=1.5pt, fill=black,circle,draw]
%\tikzstyle{every node}=[circle,draw]%fill=blue!20]

\node (a) at (0,4) {};
\node (b) at (-3.464, -2) {};
\node (c) at (3.464,-2) {};
\node (ab1) at (-2.0784,1.2) {};
\node (ab2) at (-2.598, 1.5) {};
\node (ac1) at (2.0784, 1.2) {};
\node (ac2) at (2.598,1.5) {};
\node (bc1) at (0,-2.4) {};
\node (bc2) at (0,-3) {};

%\tikzstyle{every node}=[fill=red!20]
\foreach \from/\to in {a/ab1,a/ab2,a/ac1,a/ac2,b/ab1,b/ab2,b/bc1,b/bc2,c/ac1,c/ac2,c/bc1,c/bc2}
\draw [-] (\from) -- (\to);
\end{tikzpicture}

\caption{A widespread graph that is not a forest of chandeliers} \label{fig:4}
\end{figure}
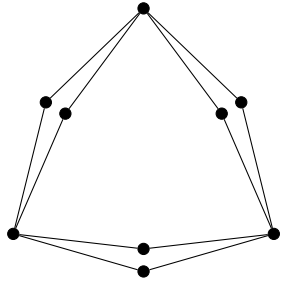

So we would like to determine which graphs are widespread. In view of the conjecture \ref{chandconj},
we expected the answer to be, again, forests of chandeliers, but this is false. 
In~\cite{bananatrees}, we found a multigraph that is not a forest
of chandeliers, that is widespread: the multigraph obtained from a triangle by adding many parallel edges 
between two of the three pairs of vertices,
and adding one parallel edge (so, two total) between the last pair. (This can be converted to a graph with the same properties by 
subdividing edges.) This does not disprove the conjecture \ref{chandconj}, but it does mean we have
no idea which graphs are widespread. Perhaps they all are?  As discussed in \cite{bananatrees}, the case of 2-controlled graphs is critical here; indeed, the conjecture that all graphs are widespread is equivalent to the following.

\begin{thm}\label{2controlconj}
{\bf Conjecture: }For every graph $H$, the ideal of all $H$-subdivision-free graphs is $2$-controlled.
\end{thm}

\section{Graphs with geometric representations}\label{sec:geometric}

A particularly interesting class of problems arises when we consider graphs that arise from geometric constructions.
Many of these are {\em intersection graphs}: given a collection $\mathcal F$ of sets, the intersection graph 
$I(\mathcal F)$ has vertex set $\mathcal F$, and distinct $X,Y\in\mathcal F$ are adjacent whenever $X\cap Y$ is 
nonempty.  For example, if we begin with a collection of intervals in the real line, then the corresponding 
intersection graph is an {\em interval graph}.  These are well known to be perfect, and therefore $\chi$-bounded.  
More generally, a {\em $d$-dimensional box graph} is an intersection graph of a collection of boxes (i.e.~products 
of intervals) in ${\mathbb R}^d$.  So 1-dimensional box graphs are interval graphs.
Asplund and Gr\"unbaum~\cite{asplund}  showed that the ideal of 2-dimensional box graphs is  $\chi$-bounded.  
Surprisingly, the ideal of 3-dimensional box graphs is not $\chi$-bounded: this follows from the construction of 
Burling discussed in section~\ref{sec:examples}, which can be realised as a 3-dimensional box graph.  
See Gy\'arf\'as and Lehel~\cite{gyarfaslehel} and Gy\'arf\'as~\cite{gyarfascorr} for related results.

Burling's construction has repeatedly proved useful. For example, another generalization of interval graphs is the 
ideal of intersection graphs of line segments in the plane.  In the 1970s, Erd\H os asked whether this ideal 
contains triangle-free graphs with arbitrarily large chromatic number.  This was resolved only recently: the 
beautiful paper of Pawlik, Kozik, Krawczyk, Laso\'{n}, Micek, Trotter and Walczak~\cite{sevenpoles} (see 
also~\cite{sixpoles}) already mentioned in section~\ref{subdivide} shows that Burling's construction can be 
realised as an intersection graph of line segments in the plane, thus answering the question of Erd\H os in the 
positive, and showing that this ideal is not $\chi$-bounded.

An important, and more general, ideal is the ideal of string graphs (also discussed in the previous section).  
A {\em string graph} is the intersection 
graph of a collection of curves in the plane (see, for instance,~\cite{matousek}).  The ideal of string 
graphs is clearly not $\chi$-bounded, as it contains the ideal of intersection graphs of line segments in the 
plane, but various interesting subideals are.  McGuinness (see~\cite{McG1,McG2,McG3}) made a significant study of 
intersection graphs of objects in the plane.  Among other things, he showed that the ideal of 
triangle-free intersection graphs of simple curves that cross a fixed line in exactly one point have bounded 
chromatic number.  This was generalized by Suk~\cite{suk}, who showed that the ideal of intersection graphs of 
simple families of curves that intersect the $y$-axis and do not intersect any vertical line in more than one 
point is $\chi$-bounded.  A nice consequence of this is the following:
\begin{thm}
The ideal of intersection graphs of unit segments in the plane is $\chi$-bounded.  
\end{thm}
 Further progress was made by Laso\'n, Micek, Pawlik and Walczak~\cite{lason},   Rok and Walczak~\cite{rok1} and Suk and Walczak~\cite{sukw}.  A very general result was proved by Rok and Walczak~\cite{rok2}:
 
\begin{thm} 
For every integer $t \ge1$, the ideal of intersection graphs of curves each crossing a fixed curve in at least 
one and at most $t$ points is $\chi$-bounded.
\end{thm}

A particularly useful property of string graphs in our context is that if they have large chromatic number then 
some small ball has large chromatic number.  The following was proved in \cite{strings}, building on and 
simplifying ideas of McGuinness \cite{McG2}.
\begin{thm}\label{stringcontrol}
The ideal of string graphs is 2-controlled.
\end{thm}
This is an important ingredient in the proof of \ref{stringgraphs}.

There are other interesting ways to define graphs from geometric objects.  For example,
the {\em visibility graph} of a set $S$ of points in the plane has vertex set $S$, with vertices $x,y$ adjacent 
if the line segment from $x$ to $y$ contains no other points from the set.  The class of visibility graphs is not an ideal, because visibility graphs can have induced
subgraphs that are not themselves visibility graphs (deleting a vertex can make two other vertices visible to each other).
K\'ara,  P\'or and Wood~\cite{kara} showed that visibility graphs with clique number at most four have bounded chromatic number, and
conjectured that for all $\kappa$, all visibility graphs with clique number at most $\kappa$ have bounded chromatic number.
But this was disproved by Pfender~\cite{pfender}, who showed there are visibility graphs with arbitrarily large chromatic number
and with clique number six.

\section{Connections: the Erd\H{o}s-Hajnal conjecture}

In this section, we look at cliques and stable sets of graphs in ideals defined by forbidden subgraphs.
It is well known from Ramsey theory~\cite{erdosszekeres} that every graph on $n$ vertices contains a clique or
stable set of size at least $\frac12\log n$.  This is tight up to a constant factor: considering random graphs  
shows that most graphs on $n$ vertices contain no clique or stable set of
size more than $2 \log n$~\cite{erdos47}. The celebrated Erd\H os-Hajnal
conjecture asserts that $H$-free graphs (that is, graphs that do not contain an induced copy of $H$) have 
much larger cliques or stable sets. Let us say that an
ideal $\mathcal{I}$ has the {\em Erd\H os-Hajnal property}  if there is some $\epsilon > 0$ 
such that every graph $G \in \mathcal{I}$ has a clique or stable set of size at least $|G|^\epsilon$. 
The Erd\H{o}s-Hajnal conjecture~\cite{EH0, EH} is the following:
\begin{thm}\label{ehconj}
{\bf Conjecture:} For every graph $H$, the ideal of $H$-free graphs has the Erd\H os-Hajnal property.
\end{thm}
There is a natural connection here to $\chi$-boundedness: if an ideal $\mathcal{I}$ is $\chi$-bounded with a 
{\em polynomial} $\chi$-binding function $f$ then $\mathcal{I}$ satisfies the Erd\H os-Hajnal property, as every 
graph $G\in \mathcal{I}$ satisfies 
$$\alpha(G)\ge|G|/\chi(G)\ge |G|/f(\omega(G))$$ 
and so $\alpha(G)f(\omega(G))\ge|G|$.  There is no implication in the other direction.  For example, the 
ideal of triangle-free graphs has the Erd\H os-Hajnal property, but is not $\chi$-bounded.

There is an interesting property that is stronger than the Erd\H os-Hajnal property:
an ideal $\mathcal{I}$ has the {\em strong Erd\H{o}s-Hajnal property} if there exists $\vare>0$
such that for every graph $G\in \mathcal{I}$ with $|G|>1$, there exist disjoint $A,B\subseteq V(G)$ with $|A|,|B|\ge \vare|G|$
such that $A,B$ are complete or anticomplete. 
Here, we say that two disjoint sets $A,B$ are {\em complete} if every vertex in $A$ is adjacent to every 
vertex in $B$, and {\em anticomplete} if there are no edges between $A,B$.
It is not hard to show that if an ideal has the strong Erd\H{o}s-Hajnal property then it has 
the Erd\H{o}s-Hajnal property (see \cite{alontosharir,foxpach,sparselinear}).

Which ideals have the strong Erd\H{o}s-Hajnal property? Let us start with ideals defined by excluding one graph.  
The random graph of Erd\H{o}s described in section~\ref{sec:examples} has, with high probability, no pair of 
linear-sized sets that are either complete or anticomplete; so for all $g$, there is a graph not in $\mathcal{I}$ 
and 
with no cycle of length at most $g$.
Consequently, if the ideal of all $H$-free graphs has the strong Erd\H{o}s-Hajnal property, then $H$ has no cycles, and nor
does its complement; so $H$ has at most four vertices. Thus, for ideals defined by excluding one graph, the 
strong Erd\H{o}s-Hajnal property is not very interesting.

What about ideals defined by excluding a finite set of induced subgraphs?
Once again, because of Erd\H{o}s' random graph, one of the excluded subgraphs must be a forest and one
must be the complement of a forest.
In this case, an interesting result was proved by Bousquet, Lagoutte and Thomass\'e~\cite{lagoutte}:
\begin{thm}\label{pathEH}
For every path $H$,
the ideal of all graphs that contain neither $H$ nor $\overline{H}$ has the strong Erd\H{o}s-Hajnal property. 
\end{thm}
This was extended by Choromanski, Falik, Liebenau, Patel and Pilipczuk~\cite{hooks}, who proved the same for
trees $H$ formed from a path by adding a leaf adjacent to its third vertex.
Which other graphs $H$ can we take here? It is clearly necessary for one of $H$, $\overline H$ to be a forest. 
Liebenau and Pilipczuk~\cite{liebenau} conjectured that this is sufficient:
that for every forest $H$, the ideal of all $\{H,\overline{H}\}$-free graphs has the strong Erd\H{o}s-Hajnal 
property. This has recently been proved, as we discuss below.

%This is still rather far from the Gy\'arf\'as-Sumner conjecture, but it gets closer. 
A useful tool in these problems is a theorem of R\"odl~\cite{rodl}, which 
says:
\begin{thm}\label{rodl}
For every graph $H$ and all $\vare>0$ there exists $\delta>0$ such that for every $H$-free graph $G$,
there exists $X\subseteq V(G)$ with $|X|\ge \delta|G|$ such that in one of $G[X]$, $\overline{G}[X]$,
every vertex in $X$ has degree at most $\vare|X|$.
\end{thm}
The proof is a straightforward regularity lemma argument, but the result is important as it allows us to concentrate on the sparse and dense cases.
If we are excluding a pair $H,\overline H$ then the problem is invariant under complementing $G$.  So to prove 
that the ideal of all $\{H,\overline{H}\}$-free graphs has the strong Erd\H{o}s-Hajnal property,
it suffices to show the ``one-sided'' result that for some $\vare>0$, 
if $G\in \mathcal{I}$ has at least two vertices and has maximum degree less than $\vare |G|$, then 
there exist disjoint $A,B\subseteq V(G)$ with $|A|,|B|\ge \vare|G|$
such that $A,B$ are anticomplete  (``complete'' is not an outcome since the maximum degree is less than $\vare |G|$).  

In~\cite{cats}, Liebenau, Pilipczuk, Spirkl and the second author proposed a strengthening of the 
Liebenau-Pilipczuk conjecture, 
the statement~\ref{trees} below; 
and proved it for ``subdivided caterpillars'', 
that is, trees such that all their vertices of degree more than two lie in a path.
Very recently, we (with Chudnovsky and Spirkl)~\cite{trees} have proved this conjecture:
\begin{thm}\label{trees}
For every forest $H$, there exists $\vare>0$
such that for every $H$-free graph $G$ with at least two vertices, either some vertex has degree at least $\vare|G|$
or there exist disjoint $A,B\subseteq V(G)$ with $|A|,|B|\ge \vare|G|$
such that $A,B$ are anticomplete.
\end{thm}
Once again, this is best possible, because no graph $H$ that is not a forest has this property.
Using~\ref{rodl}, this implies the
Liebenau-Pilipczuk conjecture, in full; that is:
\begin{thm}\label{hkEH}
For all forests $H,K$, the ideal of all graphs that contain neither $H$ nor $\overline K$ has the strong Erd\H os-Hajnal property.  
\end{thm}
We feel this gives a little support to the Gy\'arf\'as-Sumner conjecture. It shows that something major happens 
when a forest is excluded; not what Gy\'arf\'as-Sumner predicts, but something. 
It also suggests that it might be fruitful to look at the following weakening of the Gy\'arf\'as-Sumner 
conjecture, raised by Gy\'arf\'as in~\cite{gyarfas}.
\begin{thm}\label{weakgs}
{\bf Conjecture: }For all forests $H$, the ideal of all graphs that contain neither $H$ nor $\overline H$ 
is $\chi$-bounded.
\end{thm}

Another interesting parallel with $\chi$-boundedness comes when we exclude induced subdivisions of a graph.
Recall from section~\ref{subdivide} that if $H$ is a tree or a cycle (or more generally a banana tree) then the 
ideal of graphs with no induced subdivision of $H$ is $\chi$-bounded; and that this does not hold unless $H$ is a 
forest of chandeliers.  
An intriguing result of Bonamy, Bousquet and Thomass\'e, extending \ref{pathEH}, shows:
\begin{thm}\label{holeEH}
For every cycle $C$,
the ideal of all graphs $G$ such that neither $G$ nor $\overline G$  contains an induced subdivision of $C$ has the strong Erd\H{o}s-Hajnal property. 
\end{thm}
We recently (with Chudnovsky and Spirkl)~\cite{sparselinear} proved the following very substantial extension of this result.
\begin{thm}\label{subdividedEH}
For all graphs $H,K$, the ideal of all graphs $G$ such that $G$ does not contain an induced subdivision of $H$ and $\overline G$ does not contain an induced subdivision of $K$ has the strong Erd\H os-Hajnal property.  
\end{thm}
Once again, this follows from a (stronger) ``one-sided'' result analogous to \ref{trees}:
\begin{thm}\label{moretrees}
For every graph $H$, there exists $\vare>0$
such that every graph $G$ with at least two vertices contains one of the following:
\begin{itemize}
 \item an induced subdivision of $H$; 
\item a vertex of degree at least $\vare |G|$; or
 \item disjoint anticomplete sets $A,B$ of size at least $\vare|G|$.
\end{itemize}
\end{thm}
Interestingly, the proof strategy is in part adapted from some of the $\chi$-boundedness arguments discussed in 
sections~\ref{sec:holes} and~\ref{subdivide}.  For further discussion on the connections, see~\cite{sparselinear}.

\section{Gy\'arf\'as' complementation conjecture}\label{sec:comp}

Let us mention one other result about $\chi$-boundedness. The idea of $\chi$-boundedness grew as a generalization of 
perfect graphs; and the complement of a perfect graph is perfect. But if an ideal $\mathcal{I}$ is $\chi$-bounded, 
it does not follow that the ideal of
all complements of members of $\mathcal{I}$ is $\chi$-bounded. For instance, the ideal of all graphs with stability number 
at most two is
$\chi$-bounded ($\chi(G)\le |G|\le \omega(G)+\omega(G)^2$, for all graphs $G$ in the ideal);
but the ideal of their complements is not; while
the ideal of intersection graphs of line segments in the plane is not $\chi$-bounded, but the ideal of complements 
of these graphs is (this follows from work of Pach and T\"or\H{o}csik~\cite{PaTo}).

Let $f$ be a function from the set of nonnegative integers into itself.
We recall that $f$ is a {\em $\chi$-binding function} for an ideal $\mathcal{I}$ if $\chi(G)\le f(\omega(G))$ for every
$G\in \mathcal{I}$. Let $\mathcal{I}_f$ be the maximal ideal with $\chi$-binding function $f$ (it is unique).
We say $f$ has a {\em complementary $\chi$-binding function} if the ideal of complements of members of $\mathcal{I}_f$
is $\chi$-bounded. Thus, Lov\'asz's theorem~\cite{lovaszperfect}, that the complement of a perfect graph is perfect, 
implies that the function $f(x)=x$ has a complementary $\chi$-binding function,
and we might ask which other functions do.

Gy\'arf\'as~\cite{gyarfas} proved that if $f$ has a complementary $\chi$-binding function then $\inf_{x\rightarrow \infty} f(x)/x=1$
(see~\cite{nearlyperfect} for a sharpening), and conjectured
that for all $c$ the function $f(x)=x+c$ has a complementary $\chi$-binding function. We proved this
in~\cite{gyarfascomp}:

\begin{thm}\label{gyarfascomp}
Let $c\ge 0$, and let $\mathcal{I}$ be an ideal such that
$\chi(G)\le \omega(G)+c$ for all $G\in \mathcal{I}$. Let $\mathcal{I}'$ be the ideal of
complements of members of $\mathcal{I}$. Then $\mathcal{I}'$ is $\chi$-bounded.
\end{thm}

The proof is quite nice, so here is a sketch. No graph in $\mathcal{I}'$ has $c+1$ odd holes, 
pairwise anticomplete;
for otherwise there would be a graph in $\mathcal{I}$ consisting of $c+1$ odd antiholes, pairwise complete, and 
this graph has chromatic number 
$c+1$ more than its clique number, contrary to the hypothesis. This turns out to be all we need; we prove that
every ideal of graphs without $c+1$ odd holes, pairwise anticomplete, is $\chi$-bounded. Thus, we need to show that if $G$
does not have $c+1$ odd holes, pairwise anticomplete, and has clique number at most $\kappa$,
then $\chi(G)\le f(c,\kappa)$, for some appropriate function $f(c,\kappa)$. This follows from ~\ref{multiholeparity}, but at the time 
we had not proved the latter, and the direct proof we found is worth explaining. We use induction on $c+\kappa$. The result
holds if $c=0$, 
by \ref{oddholes}, and trivially if $\kappa=0$, so we may assume that $c,\kappa>0$.

Thus we can assume that $G$ has an odd hole; let $C$ be a shortest odd hole. From the inductive hypothesis, the set of vertices
with no neighbour in $C$ has chromatic number at most $f(c-1,\kappa)$; so we just need to show that the set of vertices 
that have neighbours in $C$  has bounded chromatic 
number. If $C$ has bounded length, then this is true, since for every vertex in $C$, its set of neighbours has bounded
chromatic number, from the inductive hypothesis. So we can assume $|C|$ is as large as we like.
But $C$ is a shortest odd hole, so we can say a lot about the vertices with neighbours in it. Indeed,
the algorithm of~\cite{bergealg} was mostly concerned with these neighbours, and we can use many of the
same ideas here. We omit further details.

\section{Operations on $\chi$-bounded ideals}\label{sec:operations}

In this section we consider the closure of ideals under various operations.  We shall be interested in operations that act on a 
single graph, or a finite set of graphs, and produce a graph as output.  The {\em closure} of an ideal $\mathcal G$ under an 
operation $\phi$ is the smallest ideal that is closed under $\phi$ and contains $\mathcal G$
(or equivalently, the smallest class that contains $\mathcal G$ and is closed under $\phi$ and taking induced subgraphs).  
The main question is,  under which operations does the 
closure of a $\chi$-bounded ideal $\mathcal G$ remain $\chi$-bounded?  In many cases, this is equivalent to considering an 
ideal $\mathcal G^*$ and some way to decompose elements of $\mathcal G^*$ into smaller graphs in $\mathcal G^*$.  If we write 
$\mathcal G$ for the basic (indecomposable) graphs in $\mathcal G^*$, when does the $\chi$-boundedness of $\mathcal G$ imply the $\chi$-boundedness of $\mathcal G^*$?

For example, it is easy to see that if we take the closure of a $\chi$-bounded ideal $\mathcal G$ under the operation of 
gluing along cliques, then the resulting ideal is $\chi$-bounded (with the same $\chi$-binding function).  

A more interesting example is given by substitution.  Given graphs $G,H$ and a vertex $v\in V(G)$,  the graph obtained by 
{\em substituting $H$ for $v$} is obtained from $G$ by deleting $v$ from $G$, adding a vertex-disjoint copy of $H$, and joining 
every vertex of $H$ to every neighbour of $v$ in $G$.   Lov\'asz \cite{Lov83} showed that the ideal of perfect graphs is closed 
under substitution.  In general, $\chi$-bounded classes are not closed under substitution, but the following was shown in \cite{CPST}:

\begin{thm}
If $\mathcal G$ is a $\chi$-bounded ideal then the closure of $\mathcal G$ under substitution is also $\chi$-bounded.  Furthermore, 
if $\mathcal G$ has a polynomial $\chi$-binding function then so does its closure; and if $\mathcal G$ has an exponential 
$\chi$-binding function then so does its closure.
\end{thm}

Other operations are known to preserve $\chi$-boundedness.  For example, the closure of an ideal under the operation of gluing along 
at most $k$ vertices (where the two graphs being glued together have the same induced subgraph on the overlap) preserves 
$\chi$-boundedness (this is proved in  \cite{CPST}, or can be deduced from from earlier work of Alon, Kleitman, Saks,  Seymour, 
and Thomassen \cite{AKSST}).   
The closure of a $\chi$-bounded ideal under 1-joins is also $\chi$-bounded
(see Dvo\v r\'ak and {Kr\'al\textquoteright}~\cite{dvorak}, Bonamy and Pilipczuk~\cite{bonamy3}, and
Kim, Kwon, Oum and Sivaraman~\cite{kim}).  It would be interesting to know what happens with other graph compositions 
(see~\cite{CPST} for discussion).

A little is known about combining operations. It is known that the closure under both substitution and gluing along cliques 
preserves $\chi$-boundedness, as does the closure under gluing along sets of bounded size and cliques \cite{CPST}.  A more general 
result was proved by Penev~\cite{penev}, who shows that it is possible to add in amalgams as well.  
However the following question from~\cite{CPST} is open.
\begin{thm}
{\bf Question:}
Is the closure of a $\chi$-bounded ideal under substitution and gluing along a bounded number of vertices $\chi$-bounded?
\end{thm}

In general, given two operations that separately preserve $\chi$-boundedness of ideals when taking the closure, it is not clear 
that $\chi$-boundedness is preserved when taking the closure under both operations together.  However,  the following generalizes a problem raised in~\cite{CPST}.

\begin{thm}\label{close2}
{\bf Problem:}
Suppose that $\phi_1$ and $\phi_2$ are operations such that the closure of every $\chi$-bounded ideal under either $\phi_1$ or 
$\phi_2$ is $\chi$-bounded.  Is the closure of every $\chi$-bounded ideal under $\phi_1$ and $\phi_2$ $\chi$-bounded?
\end{thm}

Note that it is important that the closure is an ideal: if we instead consider just the smallest class closed under applying the operations and do not demand that it is closed under taking induced subgraphs then the problem has a negative answer (see \cite{CPST}, which uses this weaker form of closure).

\section{Open problems}\label{openproblems}

Here are some open questions, different from those discussed earlier. (There are yet more in~\cite{sivaraman}.)

\subsection*{Triangle-free subgraphs}

There is a fundamental conjecture, due to Louis Esperet (unpublished), although it might have been asked before.
(See~\cite{trot}, where they say ``This has been
asked several times by researchers, but we could not find a reference.'')
\begin{thm}\label{trot}
{\bf Conjecture:} For all $\kappa,n$, every graph with sufficiently large chromatic number
and clique number at most $\kappa$ has a triangle-free induced subgraph with chromatic number at least~$n$.
\end{thm}
If we relax our requirements and do not demand that the subgraph be induced, then the conjecture is true: R\"odl \cite{rodl2}  proved that, for all $n$, 
every graph with sufficiently large chromatic number contains a triangle-free subgraph with chromatic number at least $n$.  (The stronger conjecture
of Erd\H os and Hajnal  \cite{ehn}
that every graph with huge chromatic number contains a subgraph with large girth and large chromatic number is still open.)

Conjecture \ref{trot} would serve
to reduce many questions about $\chi$-boundedness to questions about triangle-free graphs, which might be much easier.
More exactly, it would imply that an ideal is $\chi$-bounded if and only if the                                                  
triangle-free graphs in it have bounded chromatic number.\footnote{Note that it is important here that we are considering an ideal: for 
example the class of graphs of form $K_3\cup H$, where $H$ is triangle-free, is not $\chi$-bounded, and it contains no triangle-free graphs; but then it
is not closed under taking induced subgraphs.  A similar trick, taking a sequence of graphs with increasing clique number and much more rapidly increasing chromatic number, can be used to construct a class of graphs where the optimal $\chi$-binding function grows at arbitrary rate.}
For instance, every triangle-free graph with chromatic number at least three has an odd hole, so \ref{oddholes}
would be immediate; and the conjecture before \ref{holeseq} would be reduced to
\ref{holeseq} itself.
This is a great question, and we do not
even know whether it holds for $n=4$ (it is true for $n=3$ by \ref{oddholes}).

\subsection*{Polynomial $\chi$-boundedness}

Another question from Esperet~\cite{esperet}: an ideal is {\em polynomially $\chi$-bounded} if it has a polynomial $\chi$-binding function 
(defined in section \ref{sec:comp}). We have no example of an ideal that is known to be $\chi$-bounded and not polynomially 
$\chi$-bounded (although there are many that are known to be $\chi$-bounded and not known to be polynomially $\chi$-bounded).
Could it be true that every $\chi$-bounded ideal is polynomially $\chi$-bounded? This would imply that every $\chi$-bounded ideal satisfies the Erd\H{o}s-Hajnal conjecture.  

Some ideals of graphs are known to be polynomially $\chi$-bounded (see for instance 
Schiermeyer and Randerath~\cite{schiermeyer}, and Bonamy and Pilipczuk~\cite{bonamy3}).
%Gajarsk\'y, Kreutzer, Ne\v set\v ril,
%Ossona de Mendez, Pilipczuk, Siebertz and Torunczyk \cite{icalp}).
But it is easy to see that there is no $d$ such that every $\chi$-bounded ideal is bounded by a polynomial of degree $d$.  
(Consider the ideal of graphs with no stable set of size $t$.)
As noted by Trotignon and Pham~\cite{trotignon} (which also contains interesting further discussion),
even the following problem is open.
\begin{thm}
 {\bf Problem.}  Is it true that, for every $t$, the ideal $\mathcal{I}_t$ of graphs with no induced $t$-vertex path is polynomially $\chi$-bounded?
\end{thm}
(This is open even for $t=5$, and its truth in that case would settle the smallest open case of the Erd\H{o}s-Hajnal conjecture.)
Gy\' arf\'as~\cite{gyarfas} shows that 
every $\chi$-binding function $f$ for $\mathcal{I}_t$ must satisfy
$$f(\omega)\ge \frac{R(\lceil t/2\rceil,\omega+1)-1}{\lceil t/2\rceil-1},$$
as a graph $G$ with no stable set of size $\lceil t/2\rceil$ does not contain an induced path on 
$t$ vertices and has chromatic number at least $|G|/\alpha(G)\ge |G|/({\lceil t/2\rceil-1})$.
By results of Spencer~\cite{spencer} on off-diagonal Ramsey numbers (see Bohman and Keevash~\cite{bohmankeevash} for a strengthening), 
$f(\omega)=\Omega((\omega/\log\omega)^{(t+1)/4})$ for fixed $t$, so any polynomial $\chi$-binding function for $\mathcal{I}_t$ must have
degree at least $(t+1)/4$.
From the other side, \ref{paths} implies that $\mathcal{I}_t$ is  
$\chi$-bounded, and that proof gives a singly-exponential $\chi$-binding function for $\mathcal{I}_t$.
The best result we know on these lines is due to Gravier, Ho\`ang and Maffray~\cite{gravier}:
\begin{thm}\label{gravierthm}
For all $t\ge 4$, if a graph $G$ has no induced $t$-vertex path, then $\chi(G)\le (t-2)^{\omega(G)-1}$.
\end{thm}

Trotignon and Pham~\cite{trotignon}  raised another intriguing question on the growth of $\chi$-binding functions.  
Given a $\chi$-bounded ideal $\mathcal G$, we can define an ``optimal'' $\chi$-binding function $f_{\mathcal G}$ by 
$$f_{\mathcal G}(k)=\max\{\chi(G):G\in\mathcal G, \omega(G)\le k\}.$$
Suppose we know the value of $f_{\mathcal G}(2)$.  Does this tell us anything about the possible values of $f_{\mathcal G}(3)$?  
If $f_{\mathcal G}(2)=2$ then $\mathcal G$ is contained in the ideal of graphs with no odd hole, and so $f_{\mathcal G}(3)$ is bounded 
(in fact, a result of Chudnovsky, Robertson, Seymour and Thomas \cite{k4free} implies that $f_{\mathcal G}(3)\le 4$).  But what 
if $f_{\mathcal G}(2)=3$?  Can $f_{\mathcal G}(3)$ be arbitrarily large? This would be answered by \ref{trot}, but might be easier
than solving \ref{trot} in general.

\subsection*{Rainbow subgraphs}

When we were discussing rainbow subgraphs before, we coloured the graph with an arbitrary number of colours; but what if we colour it
optimally? Thus, which graphs $H$ have the following property? Let 
$G$ be a graph with huge chromatic number and bounded clique number, coloured with $\chi(G)$ colours; then a rainbow copy of $H$
is present as an induced subgraph.

Again, all such graphs $H$ are forests, but more than just paths; at least all stars and all paths have this property.
Could it be that every forest has the property?

\subsection*{Hole lengths}

As discussed in section \ref{sec:holes}, the three Gy\'arf\'as conjectures on holes are now resolved; indeed, \ref{holeparitya} gives a much stronger result.  But many further questions remain.   
Let us say that a set $S$ of integers is {\em constricting} if the ideal of graphs that do not contain a hole with length in $S$ is $\chi$-bounded.  What sets are constricting?  Conjecture \ref{holeseqconj} would imply the following:

\begin{thm}\label{holeseqconj2}
{\bf Conjecture:}
Every infinite set $S$ of natural numbers with bounded gaps is constricting.
\end{thm}
Could we make do with even sparser sets of integers?  Perhaps the following wild conjecture is true:
\begin{thm}\label{constricting}
{\bf Conjecture:} A set of integers is constricting if and only if it has strictly positive lower density.
\end{thm}
Note that we cannot replace lower density by upper density here: choose a sequence $G_1,G_2,\dots$ of triangle-free graphs with 
$\chi(G_i)>i$ and with girth growing sufficiently rapidly (say $g(G_{i+1})>2^{|G_i|}$), and let $\mathcal G$ be the set of all 
induced subgraphs of the $G_i$.  Then $\mathcal G$ is clearly not $\chi$-bounded, but the set of hole lengths that do not occur 
in $\mathcal G$ has upper density 1.

We are far from proving \ref{constricting} in either direction, even for triangle-free graphs.  
At the moment, we do not even know whether the following conjecture (which contradicts \ref{constricting}!) is true.
\begin{thm}\label{densityzero}
{\bf Conjecture:}  There is a set of integers with upper density $0$ that is constricting.
\end{thm}
It would be very interesting to resolve this.  Note that we can assume that $A$ contains all integers between 3 and some large constant, so that we are considering graphs of large girth.
In the case where we do not demand that our cycles are induced, more is known:  Verstraete \cite{verstraete} showed that there is a set $A$ of
density 0 such that every graph with sufficiently large minimum degree contains (as a subgraph) a cycle whose length belongs to $A$.

Another intriguing line of inquiry is to look at whether the holes we get are short or long. Let us focus on triangle-free graphs 
(as the problems are already hard enough there).  Of course, a triangle-free graph with large chromatic number may not contain 
short holes; but if it does not, then perhaps we get holes of many different lengths.  One way to say this is the following:  

\begin{thm}\label{manyholes}{\bf Problem:}
Are there functions $f,g:\mathbb N\to \mathbb N$ such that $f(t)\to0$ and $g(t)\to\infty$ as $t\to\infty$ and the following holds?
Let $G$ be a triangle-free graph and let $C$ be the set of lengths of holes in $G$.  Then
$$\sum_{t\in C}f(t)\ge g(\chi(G)).$$
\end{thm}

Similarly, perhaps there is a version of \ref{holeseq} that takes account of the length of the holes:

\begin{thm}\label{holeseqprob2}{\bf Problem:}
Is there a function $f:\mathbb N\to \mathbb N$ such that $f(t)\to\infty$ as $t\to\infty$ and the following holds:
for every triangle-free graph $G$ of sufficiently large chromatic number, there is some $t$ such that $G$ contains holes of $f(t)$ 
consecutive lengths, each at most $t$?
\end{thm}

As a first step towards \ref{manyholes} and \ref{holeseqprob2}, surely the following variant of \ref{holeseq} must be true:

\begin{thm}\label{holeseqprob3}
{\bf Conjecture:} There is an integer $k$ such that the following holds.  For every integer $t$, every graph with chromatic number at least $k$ and sufficiently large girth contains holes of $t$ consecutive lengths.
\end{thm}

\subsection*{Excluding $K_{n,n}$}

A helpful weakening of several of these problems is to exclude $K_{n,n}$ as well as $K_n$ (or equivalently,
to just exclude $K_{n,n}$ as a subgraph rather than an induced subgraph). For instance, the following weakening
of the Gy\'arf\'as-Sumner conjecture was proved 
by Hajnal and R\"odl independently (see~\cite{gst}).
\begin{thm}\label{rodltree}
For every tree $T$ and every integer $n\ge 0$, all $T$-free graphs
not containing $K_{n,n}$ as a subgraph (not necessarily as an induced subgraph) have bounded chromatic number.
\end{thm}
There is also a beautiful theorem of K\"uhn and Osthus~\cite{kuhn}:
%that generalizes this:

\begin{thm}\label{kuhn}
For all integers $n\ge 0$,  every graph $G$ with large enough
average degree has either a $K_{n,n}$ subgraph or a subdivision
of $K_n$  as an induced subgraph.
\end{thm}
This is best possible in some sense, because $G$ itself could be a complete bipartite graph and then it contains
nothing else (see also~\cite{dvorak2}).
But if we assume that
$G$ has large chromatic number
instead of just large average degree, then we might hope to prove something stronger. Here is one possibility.
By an ``odd subdivision'' we mean a
subdivision where every edge is replaced by a path of odd length.
We hope that
\begin{thm}\label{betterkuhn}
{\bf Conjecture:} For every graph $G$, if its chromatic number is large enough, then either $G$ has a $K_{n,n}$ subgraph or
$G$ contains an odd subdivision
of $K_n$  as an induced subgraph.
\end{thm}
Using our standard approach via $r$-controlled
ideals, the problem breaks into subproblems, and we have handled all except one of them, so this might come out.

\subsection*{Colouring graphs with no long holes}
Say $G$ is {\em short-holed} if every hole in $G$ has length four.
Here is a nice question, suggested by Vaidy Sivaraman (private communication). We know that the ideal of short-holed graphs  
is $\chi$-bounded, but how well can we bound chromatic number in terms of clique number? We were able to 
prove (thanks to Sivaraman for discussions on this):
\begin{thm}\label{5holes}
If $G$ is short-holed, then $\chi(G)\le 10^{20} 2^{\omega(G)^2}$.
\end{thm}
This is no doubt still far from the truth; can we at least get a  bound that is singly-exponential? 
As far as we know, maybe $\chi(G)\le \omega(G)^2$ for all short-holed graphs (or even for all graphs with no odd hole).

On a related topic,
the Ho\`{a}ng-McDiarmid conjecture~\cite{hoangmcd} asserts the following:
\begin{thm}\label{hm}
 If $G$ has at least one edge and has no odd hole,
then its vertex set can be partitioned into two sets $A,B$ such that every maximum clique intersects both $A$ and $B$.
\end{thm}
(Note that odd cycles of length more than three do not admit such a partition.)
This is not proved even
for short-holed graphs.
\ref{hm} would imply a singly-exponential bound on $\chi$ in terms of $\omega$. 

Recently (private communication) Ho\`{a}ng has proposed a stronger conjecture,
that if $G$ has no odd hole, its vertex set can be partitioned into $\omega(G)$ sets each of which induces a 
perfect graph. This would imply that $\chi(G)\le \omega(G)^2$.

\subsection*{Cycles with chords}

We know from the results of section \ref{sec:holes} that, for every $\ell$, the ideal of graphs with no induced cycle of length at least $\ell$ is $\chi$-bounded.  What about cycles with chords?  

Trotignon and Vu\v skovi\'c~\cite{trotignonvuskovic} showed that the ideal of graphs that do not contain as an induced subgraph 
a cycle with a unique chord is $\chi$-bounded (such graphs satisfy $\chi(G)=\omega(G)$ if $\omega(G)>2$, so they are very 
close to perfect); and Trotignon and Pham~\cite{trotignon} showed that this also holds for the ideal 
of graphs that do not contain an induced subgraph that is a cycle of length at least five with a unique chord (both papers use structural 
decompositions and obtain polynomial bounds).
Aboulker and Bousquet~\cite{aboulkerbousquet} prove that, for $k=2,3$, the ideal of graphs that do not contain an induced subgraph 
that is a cycle with exactly $k$ chords is $\chi$-bounded.  More generally, they conjecture:
\begin{thm}
{\bf Conjecture:} For every $k$, the ideal of graphs that do not contain an induced subgraph that is a cycle with exactly $k$ chords is $\chi$-bounded.
\end{thm}
It seems possible that a stronger statement holds.
\begin{thm}
{\bf Conjecture:} For every $k$, the ideal of graphs that do not contain an induced cycle such that some vertex has at least $k$ neighbours on the cycle is $\chi$-bounded.
\end{thm}
If we do not require that the cycle is induced, then Aboulker~\cite{aboulker1} conjectures that the chromatic number is at most $k$
(see  Trotignon \cite{trotper}, Aboulker, Radovanovi\'c, Trotignon and Vu\v skovi\'c \cite{artv}, Bousquet and Thomass\'e, \cite{both} and Aboulker \cite{aboulker1} for results and discussion related to both conjectures).

\subsection*{Vertex-minors}

We say $H$ is a {\em vertex-minor} of $G$ if $H$ can be obtained from an induced subgraph of $G$ by (repeatedly) choosing a vertex  $v$
and replacing the subgraph induced on the neighbour set of $v$ by its complement. 
A ``circle graph'' is the intersection graph of chords of a circle. The ideal of all circle graphs is closed under taking vertex-minors,
and is $\chi$-bounded~\cite{gyarfascircle, gyarfascorr,kostochka}, and indeed  Davies and McCarty~\cite{rose} proved that it is
polynomially $\chi$-bounded, 
with a quadratic binding function. (Incidentally, a result of Kim, Kwon, Oum and Sivaraman~\cite{kim} shows that for every $n\ge3$, 
the ideal of graphs with no vertex-minor isomorphic to the cycle $C_n$ is polynomially $\chi$-bounded.)

Jim Geelen (see~\cite{dvorak})
conjectured the following (an ideal is {\em proper} if it is not the ideal of all graphs):
\begin{thm}\label{geelenconj}
{\bf Conjecture:} Every proper ideal that is closed under taking vertex-minors
is $\chi$-bounded. 
\end{thm}
Choi, Kwon, Oum and Wollan~\cite{choi} proved  
that every ideal closed under vertex-minors not containing all wheels is $\chi$-bounded.
This was extended by a recent result of Geelen, Kwon, McCarty and Wollan~\cite{geelenkwon}, proving the same for every ideal closed under vertex-minors 
that does not include all circle graphs.  
These result are both superceded
by a very recent result by James Davies, claiming that \ref{geelenconj} is true in general.

Geelen also conjectures (private communication) that for every proper ideal closed under vertex-minors, 
there are polynomial-time
algorithms to compute maximum clique and maximum stable set for its members.

\subsection*{Algorithms}

What about algorithms for all the $\chi$-boundedness theorems? For instance, we proved that graphs with no odd hole have
chromatic number bounded by a function of their clique number, but can we find such a colouring in polynomial time?
There is now a poly-time algorithm to test if a graph has an odd hole~\cite{oddholetest}, but even without the use of that,
we observed in~\cite{oddholes} that there is a poly-time algorithm that will either find an odd hole or
find a clique and a colouring, with the number of colours at most doubly-exponential in the size of the clique, 
as in the theorem.
It would be good to develop results like this for the other $\chi$-boundedness theorems.

\section*{Acknowledgements}
Several people have sent us corrections or improvements. We would particularly like to thank 
Pierre Aboulker, Pierre Charbit,
Louis Esperet, 
Patrice Ossona de Mendez, Sang-Il Oum, Bruce Reed, 
Vaidy Sivaraman, Sophie Spirkl and Nicolas Trotignon. Also, we would like to thank the referees for an exceptionally helpful set of comments.

\end{document}